# Brownian motion on disconnected sets, basic hypergeometric functions, and some continued fractions of Ramanujan

### Shankar Bhamidi[1], Steven N. Evans[*1] Ron Peled[†1] and Peter Ralph[‡1]

*University of California, Berkeley*

**Abstract:** Motivated by Lévy's characterization of Brownian motion on the line, we propose an analogue of Brownian motion that has as its state space an arbitrary closed subset of the line that is unbounded above and below: such a process will be a martingale, will have the identity function as its quadratic variation process, and will be "continuous" in the sense that its sample paths don't skip over points. We show that there is a unique such process, which turns out to be automatically a reversible Feller-Dynkin Markov process. We find its generator, which is a natural generalization of the operator $f \mapsto \frac{1}{2} f''$.

We then consider the special case where the state space is the self-similar set $\{\pm q^k : k \in \mathbb{Z}\} \cup \{0\}$ for some $q > 1$. Using the scaling properties of the process, we represent the Laplace transforms of various hitting times as certain continued fractions that appear in Ramanujan's "lost" notebook and evaluate these continued fractions in terms of basic hypergeometric functions (that is, $q$-analogues of classical hypergeometric functions). The process has 0 as a regular instantaneous point, and hence its sample paths can be decomposed into a Poisson process of excursions from 0 using the associated continuous local time. Using the reversibility of the process with respect to the natural measure on the state space, we find the entrance laws of the corresponding Itô excursion measure and the Laplace exponent of the inverse local time – both again in terms of basic hypergeometric functions. By combining these ingredients, we obtain explicit formulae for the resolvent of the process. We also compute the moments of the process in closed form. Some of our results involve $q$-analogues of classical distributions such as the Poisson distribution that have appeared elsewhere in the literature.

## 1. Introduction

Let $\mathbb{T}$ be an arbitrary closed subset of $\mathbb{R}$. There is a well-developed theory of differentiation, integration, and differential equations on $\mathbb{T}$ (sometimes refered to as

---

*Supported in part by NSF Grant DMS-0405778. Part of the research was conducted while the author was visiting the Pacific Institute for Mathematical Sciences, the Zentrum für interdisziplinäre Forschung der Universität Bielefeld, and the Mathematisches Forschungsinstitut Oberwolfach.

†Supported in part by a Liné and Michel Loève Fellowship.

‡Supported in part by a VIGRE grant awarded to the Department of Statistics at U.C. Berkeley.

[1]University of California, Department of Statistics, 367 Evans Hall #3860, Berkeley, CA 94720-3860, USA, e-mail: {shanky, evans, peledron, plr}@stat.berkeley.edu

the *time scale calculus*) that simultaneously generalizes the familiar Newtonian calculus when $\mathbb{T} = \mathbb{R}$ and the theory of difference operators and difference equations when $\mathbb{T} = \mathbb{Z}$ (as well as the somewhat less familiar theory of $q$-differences and $q$-difference equations when $\mathbb{T}$ is $\{q^k : k \in \mathbb{Z}\}$ for some $q > 1$). The time scale calculus is described in [6], where there is also discussion of the application of time scale dynamic equations to systems that evolve via a mixture of discrete and continuous mechanisms.

Our first aim in this paper is to investigate a possible analogue of Brownian motion with **state space** an arbitrary closed subset of $\mathbb{R}$. A celebrated theorem of Lévy says that Brownian motion on $\mathbb{R}$ is the unique $\mathbb{R}$-valued stochastic process $(\xi_t)_{t \in \mathbb{R}_+}$ such that:

(I) $\xi$ has continuous sample paths,

(II) $\xi$ is a martingale,

(III) $(\xi_t^2 - t)_{t \in \mathbb{R}_+}$ is a martingale.

A similar set of properties characterizes continuous time symmetric simple random walk on $\mathbb{Z}$ with unit jump rate: we just need to replace condition (I) by the analogous hypothesis that $\xi$ does not skip over points, that is, that all jumps are of size $\pm 1$. Note that for both $\mathbb{R}$ and $\mathbb{Z}$ the Markovianity of $\xi$ is not assumed and comes as a consequence of the hypotheses.

We show in Section 2 that on an arbitrary $\mathbb{T}$ that is *unbounded above and below* there exists a unique (in distribution) càdlàg process $\xi$ that satisfies conditions (II) and (III) plus the appropriate analogue of (I) or the "skip-free" property of simple random walk. Namely:

(I') for states $x < y < z$ in $\mathbb{T}$ and times $0 \le r < t < \infty$, if either $\xi_r = x$ and $\xi_t = z$ or $\xi_r = z$ and $\xi_t = x$, then $\xi_s = y$ for some time $s$ between $r$ and $t$.

Moreover, we demonstrate that this process is a reversible Feller-Dynkin Markov process with a generator that we explicitly compute. The proof of existence is via an explicit construction as a time change of standard Brownian motion. The proof of uniqueness (which was suggested to us by Pat Fitzsimmons) relies on a result of Chacon and Jamison, as extended by Walsh, that says, informally, if a stochastic process has the hitting distributions of a strong Markov process, then it is a time change of that Markov process.

As well as establishing the existence and uniqueness of the Brownian motion on $\mathbb{T}$ in Section 2, we give its generator, which is a natural analogue of the standard Brownian generator $f \mapsto \frac{1}{2} f''$. Note that a simple consequence of (II) and (III) is that $\xi$ has the same covariance structure as Brownian motion on $\mathbb{R}$, that is $\mathbb{E}^x[\xi_s \xi_t] - \mathbb{E}^x[\xi_s] \mathbb{E}^x[\xi_t] = s \wedge t$ for all $x \in \mathbb{T}$.

The assumption that the state space $\mathbb{T}$ is unbounded above and below is necessary. To see this, first note that $\mathbb{T}$ cannot be bounded above and below, because this would imply that if $\xi_0 = x$, then $\lim_{t \to \infty} \mathbb{E}[\xi_t^2 - t] = -\infty \ne x^2$, contradicting property (III). Assume now that $\mathbb{T}$ is unbounded above and bounded below with $\inf \mathbb{T} = a > -\infty$. Suppose $\xi_0 = x$. Choose $b \in \mathbb{T}$ with $x < b$. Put $T = \inf\{t \ge 0 : \xi_t \notin [a, b)\}$. Note by the right-continuity of $\xi$ and property (I') that $\xi_T = b$ on the event $\{T < \infty\}$. By properties (I') and (III), $\mathbb{E}[t \wedge T] = \mathbb{E}[\xi_{t \wedge T}^2] - x^2 \le a^2 \vee b^2 - x^2$, and so $T$ is indeed almost surely finite. By properties (I') and (II), $(\xi_{t \wedge T})_{t \in \mathbb{R}_+}$ is a bounded martingale with $\xi_{t \wedge T} = b$ for $t \ge T$ almost surely, but this leads to the contradiction $b = \lim_{t \to \infty} \mathbb{E}[\xi_{t \wedge T}] = x$. The proof that $\mathbb{T}$ cannot be bounded above and unbounded below is similar.

The process $\xi$ is constructed as a time-change of standard Brownian motion, a class of processes described in Itô and McKean [34], and that has been studied variously as "gap diffusions" [44], "quasidiffusions" [7, 46, 47, 48, 49, 50, 55], and (one-dimensional) "generalized diffusions" [65, 66, 67, 68]. The process $\xi$ that we study is a quasidiffusion, so results on quasidiffusions apply in this context – but it is a *distinguished* quasidiffusion among the many possible quasidiffusions taking values in $\mathbb{T}$. Quasidiffusions can exhibit behavior considerably different from that of $\xi$ – for instance, Feller and McKean [22] described a quasidiffusion that has all of $\mathbb{R}$ as its state space, but spends all its time in $\mathbb{Q}$ . These processes (with killing and appropriate boundary conditions) were shown by Löbus [54], extending work by Feller [15, 18, 20] to be the only Markov processes taking values in $\mathbb{R}$ whose generators are in some sense local, and satisfy a certain maximum principle. Various authors [47, 48, 51] have given beautiful spectral representations of quasidiffusions using Kreĭn's theory of strings [13, 19, 41].

It is natural to ask about further properties of the Brownian motion on $\mathbb{T}$. In the present paper we pursue this matter in a particularly nice special case, when $\mathbb{T} = \mathbb{T}_q := \{\pm q^k : k \in Z\} \cup \{0\}$ for some $q > 1$. In this case, the process $\xi$ started at $x$ has the same distribution as the process $(\frac{1}{q^k}\xi_{q^{2k}t})_{t \in \mathbb{R}_+}$ when $\xi$ is started at $q^k x$ for $k \in \mathbb{Z}$. This Brownian-like scaling property enables us to compute explicitly the Laplace transforms of hitting times and the resolvent of $\xi$ in terms of certain continued fractions that appear in the "lost" notebook of Ramanujan. We can, in turn, evaluate these continued fractions in terms of basic hypergeometric functions (where, for the sake of the uninitiated reader, we stress that "basic" means that such functions are the analogues of the classical hypergeometric functions to some "base" – that is, the $q$-series analogue of those functions). We recall that, in general, a $q$-analogue of a mathematical construct is a family of constructs parameterized by $q$ such that each generalizes the known construct and reduces in some sense to the known construct in the limit "$q \to 1$". This notion ranges from the very simple, such as $(q^n - 1)/(q - 1)$ being the $q$-analogue of the positive integer $n$, through to the very deep, such as certain quantum groups (which are not actually groups in the usual sense) being the "$q$-deformations" of appropriate classical groups [9, 35, 40].

For a very readable introduction to $q$-calculus see [36], and for its relation with $q$-series, see the tutorial [45], or the more extensive books [30] or [2]. What we need for our purposes is given in Section 11.

The interplay between $q$-calculus (that is, $q$-difference operators, $q$-integration, and $q$-difference equations), $q$-series (particularly basic hypergeometric functions), and probability has been explored in a number of settings both theoretical and applied. The recent paper [4] studies the connection between $q$-calculus and the exponential functional of a Poisson process

$$I_q := \int_0^\infty q^{N_t}\, dt, \quad q < 1,$$

where $N_t$ is the simple homogeneous Poisson counting process on the real line. A purely analytic treatment of the distribution of $I_q$ using $q$-calculus is given in [3]. It is interesting to note that the same functional seems to have arisen in a number of applied probability settings as well, for example, in genetics [10] and in transmission control protocols on communication networks [11]. In [42] the Euler and Heine distributions, $q$-analogues of the Poisson distribution, are studied: distributional properties are derived and some statistical applications (such as fitting these distributions to data) are explored. These analogues have arisen in contexts as varied as prior distributions for stopping time strategies when drilling for oil and studies of parasite distributions, see the references in [42]. The $q$-analogue of the Pascal distribution has also been studied in the applied context, see [43]. The properties of $q$-analogues of various classical discrete distributions are also surveyed in [52]. Both $I_q$ and the Euler distribution appear in Section 6, where they come together to form the distribution of a hitting time.

Probabilistic methods have also been used to derive various results from $q$-calculus. A number of identities (including the $q$-binomial theorem and two of Euler's fundamental partition identities) are derived in [59] by considering processes involving Bernoulli trials with variable success probabilities. Several other identities (for example, product expansions of $q$-hypergeometric functions and the Rogers–Ramanujan identities) are obtained in [58] using extensions of Blomqvist's absorption process. Some properties of $q$-random mappings are explored in [57]: in particular, the limiting probability that a $q$-random mapping does not have a fixed point is expressed via a $q$-analogue of the exponential function. Connections between $q$-series and random matrices over a finite field (resp. over a local field other than $\mathbb{R}$ or $\mathbb{C}$) are investigated in [27, 28, 29] (resp. [1, 14]).

## 2. Brownian motion on a general unbounded closed subset of $\mathbb{R}$

### 2.1. *Existence*

Let $\mathbb{T}$ be a closed subset of $\mathbb{R}$ that is unbounded above and below (that is, bilaterally unbounded). We now show existence of a Feller-Dynkin Markov process satisfying conditions (I'), (II) and (III) by explicitly constructing such a process as a time-change of Brownian motion. Let $(B_t)_{t \in \mathbb{R}_+}$ be standard Brownian motion on $\mathbb{R}$ and let $\ell_t^a$ be its local time at the point $a \in \mathbb{R}$ up to time $t \geq 0$. We choose a jointly continuous version of $\ell$ and we adopt the normalization of local time that makes $\ell$ a family of occupation densities for the Brownian motion; that is, $\int_0^t f(B_s)\, ds = \int_{\mathbb{R}} f(a)\ell_t^a\, da$ for all bounded Borel functions $f$. Equivalently, for each $a$ the process $(\ell_t^a)_{t \in \mathbb{R}_+}$ is the unique continuous non-decreasing process such that $(|B_t - a| - \ell_t^a)_{t \in \mathbb{R}_+}$ is a martingale.

We introduce the following notation from [6]. For a point $x \in \mathbb{T}$ set

$$\rho(x) := \sup\{y \in \mathbb{T} \ : \ y < x\}, \qquad \sigma(x) := \inf\{y \in \mathbb{T} \ : \ y > x\}.$$

If $\rho(x) \neq x$ say that $x$ is *left-scattered*, otherwise $x$ is *left-dense*, and similarly if $\sigma(x) \neq x$ say that $x$ is *right scattered*, otherwise $x$ is *right-dense*. Denote by $\mathbb{T}_{ss}, \mathbb{T}_{sd}, \mathbb{T}_{ds}$ and $\mathbb{T}_{dd}$ the left and right scattered, left-scattered right-dense, left-dense right-scattered and left and right dense subsets of $\mathbb{T}$, respectively.

Define a Radon measure on $\mathbb{R}$ by $\mu := \mathbf{1}_{\mathbb{T}} \cdot m + \sum_{x \in (\mathbb{T} \setminus \mathbb{T}_{dd})} \frac{\sigma(x) - \rho(x)}{2} \delta_x$, where $m$ is Lebesgue measure. Observe for any $x \in \mathbb{T}$ that $\frac{\sigma(x) - \rho(x)}{2}$ is the length of the interval of points in $\mathbb{R}$ that are closer to $x$ than to any other point of $\mathbb{T}$. Thus $\mu$ is the push-forward of $m$ by the $m$-a.e. well-defined map that takes a point in $\mathbb{R}$ to the nearest point of $\mathbb{T}$. Note that the support of $\mu$ is all of $\mathbb{T}$. Define the continuous additive functional

$$A_u^\mu := \int_{\mathbb{R}} \ell_u^a\, \mu(da)$$

and let $\theta_t^\mu$ be its right continuous inverse, that is,

$$\theta_t^\mu := \inf\{u \ : \ A_u^\mu > t\}.$$

By the *time change of $B_t$ with respect to the measure* $\mu$ ([34], §5, or [60], III.21) we mean the process $\xi_t := B_{\theta_t^\mu}$. It is easily seen that $\xi$ has $\mathbb{T}$ as its state space, and if $B_0 = x \in \mathbb{T}$ then $\xi_0 = x$ also. Moreover, it is not hard to show that $\xi$ is a Feller-Dynkin Markov process on $\mathbb{T}$.

We will need the generator of $\xi$. For that purpose, we introduce the following notation. Write

$$C_0(\mathbb{T}) := \{f : \mathbb{T} \to \mathbb{R} \ : \ f \text{ is continuous on } \mathbb{T} \text{ and tends to } 0 \text{ at infinity}\}.$$

Define a linear operator $\mathcal{G}$ on $C_0(\mathbb{T})$ as follows. For $x \in \mathbb{T}$, set $y_{x,r} := \rho(x - r)$ and $z_{x,r} := \sigma(x + r)$. Put

$$(\mathcal{G}f)(x)$$
$$:= \lim_{r\downarrow 0} \left( \frac{z_{x,r} - x}{z_{x,r} - y_{x,r}} f(y_{x,r}) + \frac{x - y_{x,r}}{z_{x,r} - y_{x,r}} f(z_{x,r}) - f(x) \right) \Big/ \left( (x - y_{x,r})(z_{x,r} - x) \right)$$
$$= \lim_{r\downarrow 0} \left( \frac{f(y_{x,r})}{(x - y_{x,r})(z_{x,r} - y_{x,r})} - \frac{f(x)}{(x - y_{x,r})(z_{x,r} - x)} + \frac{f(z_{x,r})}{(z_{x,r} - x)(z_{x,r} - y_{x,r})} \right)$$
$$= \lim_{r\downarrow 0} \left( \frac{f(y_{x,r}) - f(x)}{(x - y_{x,r})(z_{x,r} - y_{x,r})} + \frac{f(z_{x,r}) - f(x)}{(z_{x,r} - x)(z_{x,r} - y_{x,r})} \right)$$

on the domain $\mathrm{Dom}(\mathcal{G})$ consisting of those functions $f \in C_0(\mathbb{T})$ for which the limits exist for all $x \in \mathbb{T}$ and define a function in $C_0(\mathbb{T})$.

Note that $\mathcal{G}$ is a natural analogue of the standard Brownian generator $f \mapsto \frac{1}{2} f''$ and coincides with this latter operator when $\mathbb{T} = \mathbb{R}$. Note also that if $f$ is the restriction to $\mathbb{T}$ of a function that is in $C_0^2(\mathbb{R})$, then $f \in \mathrm{Dom}(\mathcal{G})$ and

$$(\mathcal{G}f)(x) = \begin{cases} \frac{f(\rho(x))}{(x-\rho(x))(\sigma(x)-\rho(x))} - \frac{f(x)}{(x-\rho(x))(\sigma(x)-x)} + \frac{f(\sigma(x))}{(\sigma(x)-x)(\sigma(x)-\rho(x))}, & x \in \mathbb{T}_{ss}, \\ \frac{f(\rho(x))-f(x)}{(x-\rho(x))^2} + \frac{f'(x)}{x-\rho(x)}, & x \in \mathbb{T}_{sd}, \\ -\frac{f'(x)}{\sigma(x)-x} + \frac{f(\sigma(x))-f(x)}{(\sigma(x)-x)^2}, & x \in \mathbb{T}_{ds}, \\ \frac{1}{2} f''(x), & x \in \mathbb{T}_{dd}. \end{cases}$$

**Proposition 2.1.** *The time change $\xi$ of standard Brownian motion $B$ with respect to the measure $\mu$ is a Feller-Dynkin Markov process on $\mathbb{T}$ that satisfies conditions (I'), (II) and (III). The generator of $\xi$ is $(\mathcal{G}, \mathrm{Dom}(\mathcal{G}))$.*

*Proof.* We have already noted that $\xi$ is a Feller-Dynkin Markov process. Given $x \in \mathbb{T}$, write $\mathbb{P}^x$ for the distribution of $\xi$ for the initial condition $\xi_0 = x$, and denote the corresponding expectation by $\mathbb{E}^x$. Under any $\mathbb{P}^x$ the property (I') is clear from the fact that the support of $\mu$ is all of $\mathbb{T}$. Before establishing properties (II) and (III) under any $\mathbb{P}^x$, we first show that the generator of $\xi$ is $(\mathcal{G}, \mathrm{Dom}(\mathcal{G}))$.

Write $(\mathcal{H}, \mathrm{Dom}(\mathcal{H}))$ for the generator of $\xi$. We begin by showing that $(\mathcal{H}, \mathrm{Dom}(\mathcal{H})) = (\mathcal{G}, \mathrm{Dom}(\mathcal{G}))$. For $x \in \mathbb{T}$ and $r > 0$, set $T_{x,r} := \inf\{t \ : \ d(\xi_t, x) > r\}$. By Dynkin's characteristic operator theorem [60], III, 12.2, $f \in \mathrm{Dom}(\mathcal{H})$, if and only if

$$(2.1) \qquad \lim_{r\downarrow 0} \frac{\mathbb{E}^x[f(\xi_{T_{x,r}})] - f(x)}{\mathbb{E}^x[T_{x,r}]}$$

exists at every $x \in \mathbb{T}$ and defines a function in $C_0(\mathbb{T})$, in which case this function is $\mathcal{H}f$.

Set $y_{x,r} := \rho(x - r)$ and $z_{x,r} := \sigma(x + r)$. Because the support of $\mu$ is all of $\mathbb{T}$,

$$\theta^\mu_{T_{x,r}} = \inf\{t \in \mathbb{R}_+ \ : \ B_t \in \{y_{x,r}, z_{x,r}\}\} =: U_{x,r}.$$

Thus

$$\mathbb{P}^x\{\xi_{T_{x,r}} = y_{x,r}\} = \frac{z_{x,r} - x}{z_{x,r} - y_{x,r}}$$

and

$$\mathbb{P}^x\{\xi_{T_{x,r}} = z_{x,r}\} = \frac{x - y_{x,r}}{z_{x,r} - y_{x,r}}.$$

Consequently,

$$\mathbb{E}^x[f(\xi_{T_{x,r}})] = \frac{z_{x,r} - x}{z_{x,r} - y_{x,r}} f(y_{x,r}) + \frac{x - y_{x,r}}{z_{x,r} - y_{x,r}} f(z_{x,r}).$$

Hence it is enough to show for all $x \in \mathbb{T}$ and $r > 0$ that

$$(2.2) \qquad \mathbb{E}^x[T_{x,r}] = (x - y_{x,r})(z_{x,r} - x).$$

Now

$$T_{x,r} = \int_{\mathbb{R}} \ell^a_{U_{x,r}} \, \mu(da) = \int_{(y_{x,r}, z_{x,r})} \ell^a_{U_{x,r}} \, \mu(da), \quad \mathbb{P}^x\text{-a.s.,}$$

and in particular,

$$\begin{aligned}
(2.3) \qquad \mathbb{E}^x[T_{x,r}] &= \int_{(y_{x,r}, z_{x,r})} \mathbb{E}^x[\ell^a_{U_{x,r}}] \, \mu(da) \\
&= \int_{(y_{x,r}, x]} \frac{2(a - y_{x,r})(z_{x,r} - x)}{z_{x,r} - y_{x,r}} \, \mu(da) \\
&\quad + \int_{(x, z_{x,r})} \frac{2(x - y_{x,r})(z_{x,r} - a)}{z_{x,r} - y_{x,r}} \, \mu(da) \\
&= \frac{2}{z_{x,r} - y_{x,r}} \bigg( (z_{x,r} - x) \int_{(y_{x,r}, x]} (a - y_{x,r}) \, \mu(da) \\
&\quad + (x - y_{x,r}) \int_{(x, z_{x,r})} (z_{x,r} - a) \, \mu(da) \bigg).
\end{aligned}$$

But, as we now show, for any points $u, v \in \mathbb{T}$, $u < v$,

$$(2.4) \qquad \int_{(u,v)} \mu(da) = v - u - \frac{\sigma(u) - u}{2} - \frac{v - \rho(v)}{2},$$

$$(2.5) \qquad \int_{(u,v)} a \, \mu(da) = \frac{v^2}{2} - \frac{u^2}{2} - u \frac{\sigma(u) - u}{2} - v \frac{v - \rho(v)}{2}$$

(note the similarity to Lebesgue integration up to boundary effects). Substituting this into (2.3) gives (2.2) after some algebra.

Let us prove the identities (2.4) and (2.5). For simplicity, we prove them in the special case when $\mathbb{T}_{ds} \cap (u,v) = \emptyset$. The proof of the general case is similar. Fix $u, v \in \mathbb{T}$, $u < v$. Since $\mathbb{T}$ is closed, we can write $(u, \rho(v)) \setminus \mathbb{T}$ as a countable union of disjoint (non-empty) open intervals $\{(a_n, b_n) \ : \ n \in \mathbb{N}\}$. We note that for any such interval $(a_n, b_n] \subseteq (u,v)$ and

$$\int_{(a_n, b_n]} \mu(da) = \frac{\sigma(b_n) - a_n}{2} = b_n - a_n - \frac{\sigma(a_n) - a_n}{2} + \frac{\sigma(b_n) - b_n}{2}.$$

Summing up over all these intervals the boundary effects cancel telescopically (since $\mathbb{T}_{ds} \cap (u, v) = \emptyset$) and, since $\sigma(a_n) = b_n$, we get

$$\sum_n \int_{(a_n, b_n]} \mu(da) = \sum_n (b_n - a_n) - \frac{\sigma(u) - u}{2} - \frac{v - \rho(v)}{2}$$
$$= \int_{\bigcup_n (a_n, b_n]} dm - \frac{\sigma(u) - u}{2} - \frac{v - \rho(v)}{2},$$

where again $m$ is Lebesgue measure. Identity (2.4) now follows since $(u, v) = \bigcup_n (a_n, b_n] \cup ((u, v) \cap \mathbb{T}_{dd})$ and, by the definition of $\mu$,

$$\int_{(u,v)} \mu(da) = \int_{(u,v) \cap \mathbb{T}_{dd}} \mu(da) + \int_{\bigcup_n (a_n, b_n]} \mu(da)$$
$$= \int_{(u,v)} dm - \frac{\sigma(u) - u}{2} - \frac{v - \rho(v)}{2}.$$

To prove identity (2.5), we note similarly that for any $n \in \mathbb{N}$

$$\int_{(a_n, b_n]} a\, \mu(da) = b_n \frac{\sigma(b_n) - a_n}{2} = \frac{b_n^2}{2} - \frac{a_n^2}{2} - a_n \frac{\sigma(a_n) - a_n}{2} + b_n \frac{\sigma(b_n) - b_n}{2},$$

so that again

$$\sum_n \int_{(a_n, b_n]} a\, \mu(da) = \sum_n \left( \frac{b_n^2}{2} - \frac{a_n^2}{2} \right) - u \frac{\sigma(u) - u}{2} - v \frac{v - \rho(v)}{2}$$
$$= \int_{\bigcup_n (a_n, b_n]} a\, dm - u \frac{\sigma(u) - u}{2} - v \frac{v - \rho(v)}{2},$$

and (2.5) follows since

$$\int_{(u,v)} a\, \mu(da) = \int_{(u,v) \cap \mathbb{T}_{dd}} a\, \mu(da) + \int_{\bigcup_n (a_n, b_n]} a\, \mu(da)$$
$$= \int_{(u,v)} a\, dm - u \frac{\sigma(u) - u}{2} - v \frac{v - \rho(v)}{2}.$$

By the Markov property of $\xi$, in order to show (II) and (III) it suffices to show that $\mathbb{E}^x[\xi_t] = x$ and $\mathbb{E}^x[\xi_t^2] = x^2 + t$ for all $x \in \mathbb{T}$ and $t \in \mathbb{R}_+$. By Dynkin's formula [60], III.10, for any $f \in \text{Dom}(\mathcal{G})$

$$(2.6) \qquad\qquad M_t := f(\xi_t) - \int_0^t (\mathcal{G}f)(\xi_s) ds$$

is a martingale (for each starting point). Note that if we formally apply the expression for $\mathcal{G}f$ to $f(x) = x$ (resp. $f(x) = x^2$), then we get $\mathcal{G}f(x) = 0$ (resp. $\mathcal{G}f(x) = 1$), and this would give properties (II) and (III) if $x \mapsto x$ and $x \mapsto x^2$ belonged to the domain of $\mathcal{G}$. Unfortunately, this is not the case, so we must resort to an approximation argument.

Fix $x \in \mathbb{T}$. Given any $r > 0$, for $R > r$ sufficiently large we have $[\rho(x - r), \sigma(x + r)] \subset (\rho(x - R), \sigma(x + R))$. For any such pair $r, R$, there are functions $g, h \in \text{Dom}(\mathcal{G})$ such that $g(w) = w$ and $h(w) = w^2$ for $w \in [\rho(x - R), \sigma(x + R)]$, and hence $\mathcal{G}g(w) = 0$ and $\mathcal{G}h(w) = 1$ for $w \in [\rho(x-r), \sigma(x+r)]$. It follows that $(\xi_{t \wedge T_{x,r}})_{t \in \mathbb{R}_+}$ and $(\xi_{t \wedge T_{x,r}}^2 - t \wedge T_{x,r})_{t \in \mathbb{R}_+}$ are both martingales under $\mathbb{P}^x$.

Hence, if $0 < r' < r''$, then

$$\mathbb{E}^x[(\xi_{t \wedge T_{x,r''}} - \xi_{t \wedge T_{x,r'}})^2] = \mathbb{E}^x[\xi_{t \wedge T_{x,r''}}^2] - \mathbb{E}^x[\xi_{t \wedge T_{x,r'}}^2]$$
$$= \mathbb{E}^x[t \wedge T_{x,r''}] - \mathbb{E}^x[t \wedge T_{x,r'}].$$

Thus $\xi_{t \wedge T_{x,r}}$ converges to $\xi_t$ in $L^2(\mathbb{P}^x)$ as $r \to \infty$, and so

$$\mathbb{E}^x[\xi_t] = \lim_{r \to \infty} \mathbb{E}^x[\xi_{t \wedge T_{x,r}}] = x$$

and

$$\mathbb{E}^x[\xi_t^2] = \lim_{r \to \infty} \mathbb{E}^x[\xi_{t \wedge T_{x,r}}^2] = x^2 + t,$$

as required. □

**Remark.** It is more standard, but we believe less natural (and equivalent) to instead regard the generator an operator on continuous functions of $\mathbb{R}$ that are linear outside the support of $\mu$. Such a generator has the natural interpretation $\frac{1}{2} \frac{d}{d\mu} \frac{d}{dx}$. These operators appear in relation to diffusion processes in Itô and McKean [34], §§5.1-5.3, and were studied further by Feller [15, 16, 17, 18, 19, 20, 21] (where the support of $\mu$ is connected), Löbus [54, 55], and Freiberg [24, 25, 26] (where $\mu$ is atomless).

## 2.2. Uniqueness

We next establish a uniqueness result that complements the existence result of Proposition 2.1.

We will apply the following result, which is a slight variant of Corollary 3.5 of [64] extending results of [8]. We make the assumption that the Markov process $X$ is a right process and that the process $Y$ is defined on a space satisfying the usual conditions to avoid listing Walsh's assumptions. We also state the result in terms of bounded rather than finite stopping times, but this is readily seen to be sufficient. The result says, roughly speaking, that if a process has the same state-dependent hitting distributions as some strong Markov process, then the process is a time-change of that Markov process. (For example, a consequence of the result is the celebrated result of Dubins and Schwarz that any continuous martingale is a time change of Brownian motion, from which Lévy's characterization of Brownian motion that we mentioned in the Introduction is an immediate corollary.)

**Theorem 2.1.** *Let $X = (\Omega, \mathcal{F}, \mathcal{F}_t, X_t, \theta_t, \mathbf{P}^x)$ be a Borel right process with Lusin state space $E$. Assume that the paths of $X$ are càdlàg and that $X$ has no traps or holding points. Let $Y$ be a càdlàg process with state space $E$ that is defined on a complete probability space $(\Sigma, \mathcal{G}, \mathbf{Q})$ equipped with a filtration $(\mathcal{G}_t)_{t \in \mathbb{R}_+}$ satisfying the usual conditions. Assume $Y_0 = x_0$ for some $x_0 \in E$ and that almost surely the sample paths of $Y$ are not constant over any time interval. Given a Borel set $B \subseteq E$, put $S_B := \inf\{t \geq 0 : X_t \in B\}$ and define the corresponding hitting kernel by $\pi_B(x, A) := \mathbf{P}^x\{X_{S_B} \in A\}$ for $x \in E$ and $A \in E$ Borel. Given a bounded $(\mathcal{G}_t)_{t \in \mathbb{R}_+}$-stopping time $T$, put $\tau = \inf\{t \geq T : Y_t \in B\}$. Suppose that*

$$\mathbf{Q}\{Y_\tau \in A \mid \mathcal{G}_T\} = \pi_B(Y_T, A)$$

*for all bounded* $(\mathcal{G}_t)_{t \in \mathbb{R}_+}$*-stopping times* $T$ *and all Borel sets* $A$ *and* $B$. *Then there exists a perfect continuous additive functional for* $X$ *with continuous inverse* $(T_t)_{t \in \mathbb{R}_+}$ *such that* $(Y_{T_t})_{t \in \mathbb{R}_+}$ *has the same distribution as* $(X_t)_{t \in \mathbb{R}_+}$ *under* $\mathbf{P}_0^x$.

**Proposition 2.2.** *Let* $\zeta$ *be a càdlàg* $\mathbb{T}$*-valued process such that* $\zeta_0 = z \in \mathbb{T}$. *Suppose that* $\zeta$ *satisfies the counterparts of properties (I'), (II) and (III) with* $\xi$ *replaced by* $\zeta$. *Then* $\zeta$ *possesses the same distribution as the particular Feller-Dynkin process* $\xi$ *of Proposition 2.1 has under* $\mathbb{P}^z$.

*Proof.* We wish to apply Theorem 2.1. Unfortunately, the process $\xi$ has holding points unless $\mathbb{T} = \mathbb{R}$. We adapt an artifice presented in Remark 1 after Theorem 3.4 in [64] to circumvent this difficulty.

Without loss of generality, we may suppose that $\zeta$ is defined on a complete probability space $(\Sigma, \mathcal{G}, \mathbf{Q})$, that this probability space is equipped with a filtration $(\mathcal{G}_t)_{t \in \mathbb{R}_+}$ satisfying the usual conditions, and that $(\zeta_t)_{t \in \mathbb{R}_+}$ and $(\zeta_t^2 - t)_{t \in \mathbb{R}_+}$ are both martingales with respect to $(\mathcal{G}_t)_{t \in \mathbb{R}_+}$. We will use $\mathbf{Q}[\cdot]$ to denote expectation with respect to the probability measure $\mathbf{Q}$.

We first show that the sample paths of $\zeta$ do not get trapped forever in any state. Given $a, b \in \mathbb{T}$ with $a < x < b$, put $R = \inf\{t \geq 0 : \zeta_t \notin (a, b)\}$. By the counterpart of property (I'), $\zeta_{t \wedge R} \in [a, b]$ and hence, by the counterpart of property (III), $\mathbf{Q}[t \wedge R] = \mathbf{Q}[\zeta_{t \wedge R}^2] - x^2 \leq a^2 \vee b^2 - x^2$. Thus $\mathbf{Q}[R] < \infty$ and, in particular, $\mathbf{Q}\{R < \infty\} = 1$. Since this is true for all $a$ and $b$, it follows that almost surely there does not exist a time $s \in \mathbb{R}_+$ and a state $y \in \mathbb{T}$ such that $\zeta_t = y$ for all $t \geq s$.

Let $S$ be a finite $(\mathcal{G}_t)_{t \in \mathbb{R}_+}$ stopping time, and put $T := \inf\{t > S : \zeta_t \neq \zeta_S\}$. It follows from the above that $T < \infty$ almost surely. Moreover, by the counterparts of properties (I') and (II) for $\zeta$ and the right-continuity of paths, $\zeta_T \in \{\rho(\zeta_S), \sigma(\zeta_S)\}$ almost surely with

$$\mathbf{Q}\{\zeta_T = \rho(\zeta_S) \,|\, \mathcal{G}_S\} = \frac{\sigma(\zeta_S) - \zeta_S}{\sigma(\zeta_S) - \rho(\zeta_S)}$$

and

$$\mathbf{Q}\{\zeta_T = \sigma(\zeta_S) \,|\, \mathcal{G}_S\} = \frac{\zeta_S - \rho(\zeta_S)}{\sigma(\zeta_S) - \rho(\zeta_S)}$$

on the event $\{\zeta_S \in \mathbb{T} \setminus \mathbb{T}_{dd}\}$. Thus $\zeta_T = \zeta_S$ almost surely on the event $\{\zeta_S \in \mathbb{T} \setminus \mathbb{T}_{ss}\}$ and hence, by the counterpart of property (III), $S = T$ almost surely on the event $\{\zeta_S \in \mathbb{T} \setminus \mathbb{T}_{ss}\}$. On the other hand, it is certainly the case that $S < T$ almost surely on the event $\{\zeta_S \in \mathbb{T}_{ss}\}$.

We next claim that, conditional on $\mathcal{G}_S$, the random variable $T - S$ is exponentially distributed with expectation $(\zeta_S - \rho(\zeta_S))(\sigma(\zeta_S) - \zeta_S)$ (where the exponential distribution with expectation 0 is of course just the point mass at 0). This must be so, of course, if $\zeta$ has the same distribution as $\xi$, and it is the key to adapting Theorem 2.1 to our setting in which the processes involved do have holding points. To see the claim, define a function $\Psi : \mathbb{T} \times \mathbb{T} \to \mathbb{R}$ by

$$\Psi(x, y) := \begin{cases} \frac{(y - \rho(x))(\sigma(x) - y)}{(x - \rho(x))(\sigma(x) - x)}, & x \in \mathbb{T}_{ss}, \\ 0, & x \in \mathbb{T} \setminus \mathbb{T}_{ss}. \end{cases}$$

Note that for each fixed $x$ the function $\Psi(x, \cdot)$ is quadratic. It follows from counterparts of properties (II) and (III) for $\zeta$ that the process

$$M_t := \mathbf{1}\{\zeta_S \in \mathbb{T}_{ss}\} \left[ \Psi(\zeta_S, \zeta_{(S+t) \wedge T}) + \frac{t \wedge (T - S)}{(\zeta_S - \rho(\zeta_S))(\sigma(\zeta_S) - \zeta_S)} \right], \quad t \in \mathbb{R}_+,$$

is a martingale with respect to the filtration $(\mathcal{G}_{S+t})_{t \in \mathbb{R}_+}$. Note that

$$M_t = \mathbf{1}\{\zeta_S \in \mathbb{T}_{ss}\} \left[ \mathbf{1}\{T - S > t\} + \int_0^t \frac{\mathbf{1}\{T - S > u\}}{(\zeta_S - \rho(\zeta_S))(\sigma(\zeta_S) - \zeta_S)} \, du \right].$$

Hence

$$\mathbf{1}\{\zeta_S \in \mathbb{T}_{ss}\} \left( \mathbf{Q}\{T - S > t \mid \mathcal{G}_S\} - 1 \right) = -\int_0^t \frac{\mathbf{1}\{\zeta_S \in \mathbb{T}_{ss}\} \mathbf{Q}\{T - S > u \mid \mathcal{G}_S\}}{(\zeta_S - \rho(\zeta_S))(\sigma(\zeta_S) - \zeta_S)} \, du,$$

and so

$$\mathbf{1}\{\zeta_S \in \mathbb{T}_{ss}\} \mathbf{Q}\{T - S > t \mid \mathcal{G}_S\} = \mathbf{1}\{\zeta_S \in \mathbb{T}_{ss}\} \exp\left( -\frac{t}{(\zeta_S - \rho(\zeta_S))(\sigma(\zeta_S) - \zeta_S)} \right),$$

as claimed.

We now apply the device from [64] mentioned above to "embellish" the process $\xi$ in order to produce a Feller-Dynkin process without traps or holding points. Set $\bar{\mathbb{T}} := (\mathbb{T}_{ss} \times \mathbb{R}) \cup ((\mathbb{T} \setminus \mathbb{T}_{ss}) \times \{0\}) \subseteq \mathbb{T} \times \mathbb{R}$. Put $U := \inf\{t > 0 : \xi_t \neq \xi_0\}$ and $C_t := t - \sup\{s < t : \xi_s \neq \xi_t\}$, with the convention $\sup \emptyset = 0$. That is, $C_t$ is the "age" of $\xi$ in the current state at time $t$. There is a Feller-Dynkin process $(\bar{\xi}_t, \bar{\mathbb{P}}^{(x,u)})$ with state-space $\bar{\mathbb{T}}$ such that under $\bar{\mathbb{P}}^{(x,u)}$ the process $\bar{\xi}$ has the same distribution as the process

$$(\xi_t, u + t), \quad 0 \le t < U,$$
$$(\xi_t, C_t), \quad t \ge U,$$

under $\mathbb{P}^x$.

Fix a Borel set $B \subseteq \bar{\mathbb{T}}$, write $\bar{T}_B := \inf\{t \in \mathbb{R}_+ : \bar{\xi}_t \in B\}$ for the first hitting time of $B$ by $\bar{\xi}$, and denote by $\bar{\pi}_B$ the corresponding hitting kernel. That is,

$$\bar{\pi}_B((x,u), A) := \bar{\mathbb{P}}^{(x,u)}\{\bar{\xi}_{\bar{T}_B} \in A\}$$

for $(x, u) \in \bar{\mathbb{T}}$ and $A$ a Borel subset of $\bar{\mathbb{T}}$. It is not hard to see that $\bar{T}_B$ is finite $\bar{\mathbb{P}}^{(x,u)}$-almost surely for all $(x, u) \in \bar{\mathbb{T}}$, and hence $\bar{\pi}_B((x,u), \cdot)$ is a probability measure concentrated on the closure of $B$ for all $(x, u) \in \bar{\mathbb{T}}$.

Let $(D_t)_{t \in \mathbb{R}_+}$ be the analogue of $(C_t)_{t \in \mathbb{R}_+}$ for $\zeta$. That is, $D_t := t - \sup\{s < t : \zeta_s \neq \zeta_t\}$. Given a finite $(\mathcal{G}_t)_{t \in \mathbb{R}_+}$ stopping time $\bar{S}$, put $\bar{T} := \inf\{t \ge \bar{S} : (\zeta_t, D_t) \in B\}$. From what we have shown above, it follows by a straightforward but slightly tedious argument that if $B$ is a finite set, then

$$(2.7) \qquad \mathbf{Q}\{(\zeta_{\bar{T}}, D_{\bar{T}}) \in A \mid \mathcal{G}_S\} = \bar{\pi}_B((\zeta_{\bar{S}}, D_{\bar{S}}), A)$$

(in particular, $\bar{T}$ is finite $\mathbf{Q}$-almost surely). If $B$ is arbitrary, then taking a countable dense subset of $B$ and writing it as an increasing union of finite sets shows that (2.7) still holds.

Theorem 2.1 gives that there is a continuous increasing process $(T_t)_{t \in \mathbb{R}_+}$ such that each $T_t$ is a $(\mathcal{G}_t)_{t \in \mathbb{R}_+}$ stopping time, $T_0 = 0$, and $((\zeta_{T_t}, D_{T_t}))_{t \in \mathbb{R}_+}$ has the same distribution as $\bar{\xi}$ under $\bar{\mathbb{P}}^{(z,0)}$, (recall that $\zeta_0 = z$). In particular, $(\zeta_{T_t})_{t \in \mathbb{R}_+}$ has the same distribution as $\xi$ under $\mathbb{P}^x$. Since property (III) holds for $\xi$ and its counterpart holds for $\zeta$, we have that $(\zeta_{T_t}^2 - t)_{t \in \mathbb{R}_+}$ is a martingale and $(\zeta_{T_t}^2 - T_t)_{t \in \mathbb{R}_+}$ is a local martingale. Thus $(T_t - t)_{t \in \mathbb{R}_+}$ is a continuous local martingale with bounded variation, and hence $T_t = t$ for all $t \in \mathbb{R}_+$, as required. $\qquad \square$

We note that, by a proof similar to that of Proposition 2.2, one can show any cádlág $\mathbb{T}$-valued process with properties ($\Gamma$) and (II) is a time-change of a process with the distribution of the process $\xi$ constructed in Proposition 2.1. It may be necessary to introduce extra randomness in the time-change to convert the holding times of the process at points in $\mathbb{T}_{ss}$ into exponential random variables, and it may also be necessary to introduce extra randomness to "complete" the sample paths of the copy of $\xi$ – as the original process may "run out of steam" and not require an entire sample path of a copy of $\xi$ to produce it (the most extreme example is a process that stays constant at its starting point). This observation is the analogue of the result of Dubins and Schwarz that any continuous martingale on the line is a time-change of some Brownian motion.

### 2.3. Reversibility

Extensions of the following result will hold more generally: under suitable hypotheses, a time-change of a Markov process that is reversible under some measure will be reversible under an appropriate new measure. Since we don't know of a suitable general reference, we provide the straightforward proof in our setting where the Markov process is Brownian motion.

**Lemma 2.1.** *The process $\xi$ of Proposition 2.1 is reversible with respect to the measure $\mu$. In particular, $\mu$ is a stationary measure for $\xi$.*

*Proof.* We have to show for all $\lambda > 0$ and all non-negative Borel functions $f$ and $g$ that

$$\int_{\mathbb{T}} f(x)\,\mathbb{E}^x\left[\int_0^\infty e^{-\lambda t} g(\xi_t)\,dt\right]\,\mu(dx) = \int_{\mathbb{T}} g(x)\,\mathbb{E}^x\left[\int_0^\infty e^{-\lambda t} f(\xi_t)\,dt\right]\,\mu(dx).$$

Now recalling $A^\mu$ is the inverse of $\theta^\mu$,

$$\int_{\mathbb{T}} f(x)\,\mathbb{E}^x\left[\int_0^\infty e^{-\lambda t} g(\xi_t)\,dt\right]\,\mu(dx) = \int_{\mathbb{T}} f(x)\,\mathbb{E}^x\left[\int_0^\infty e^{-\lambda t} g(B_{\theta_t^\mu})\,dt\right]\,\mu(dx)$$

$$= \int_{\mathbb{T}} f(x)\,\mathbb{E}^x\left[\int_0^\infty e^{-\lambda A_s^\mu} g(B_s)\,dA_s^\mu\right]\,\mu(dx)$$

$$= \int_{\mathbb{T}}\int_{\mathbb{T}} f(x)\,\mathbb{E}^x\left[\int_0^\infty e^{-\lambda \int_{\mathbb{T}} \ell_s^a\,\mu(da)}\,d\ell_s^y\right] g(y)\,\mu(dy)\,\mu(dx).$$

It follows from the reversibility of $B$ with respect to Lebesgue measure that for any $\gamma > 0$ and any non-negative bounded continuous functions $F$, $G$ and $H$,

$$\int_{\mathbb{R}}\int_{\mathbb{R}} F(x)\,\mathbb{E}^x\left[\int_0^\infty e^{-(\gamma s + \lambda \int_{\mathbb{T}} \ell_s^a H(a)\,m(da))}\,d\ell_s^y\right] G(y)\,m(dy)\,m(dx)$$

$$= \int_{\mathbb{R}} F(x)\,\mathbb{E}^x\left[\int_0^\infty e^{-(\gamma s + \lambda \int_0^s H(B_u)\,du)} G(B_s)\,ds\right]\,m(dx)$$

$$= \int_{\mathbb{R}} G(y)\,\mathbb{E}^y\left[\int_0^\infty e^{-(\gamma s + \lambda \int_0^s H(B_u)\,du)} F(B_s)\,ds\right]\,m(dy)$$

$$= \int_{\mathbb{R}}\int_{\mathbb{R}} G(y)\,\mathbb{E}^y\left[\int_0^\infty e^{-(\gamma s + \lambda \int_{\mathbb{T}} \ell_s^a H(a)\,m(da))}\,d\ell_s^x\right] F(x)\,m(dx)\,m(dy).$$

Thus (noting that each side is jointly continuous in $x$ and $y$),

$$\mathbb{E}^x\left[\int_0^\infty e^{-(\gamma s + \lambda \int_{\mathbb{T}} \ell_s^a H(a)\,m(da))}\,d\ell_s^y\right] = \mathbb{E}^y\left[\int_0^\infty e^{-(\gamma s + \lambda \int_{\mathbb{T}} \ell_s^a H(a)\,m(da))}\,d\ell_s^x\right]$$

for all $x, y \in \mathbb{R}$.

Writing $\mu$ as the vague limit of a sequence of Radon measures that have bounded density with respect $m$ and applying dominated convergence gives

$$\mathbb{E}^x \left[ \int_0^\infty e^{-(\gamma s + \lambda \int_\mathbb{T} \ell_s^a \, \mu(da))} \, d\ell_s^y \right] = \mathbb{E}^y \left[ \int_0^\infty e^{-(\gamma s + \lambda \int_\mathbb{T} \ell_s^a \, \mu(da))} \, d\ell_s^x \right]$$

for all $x, y \in \mathbb{R}$. Hence, by monotone convergence,

$$\mathbb{E}^x \left[ \int_0^\infty e^{-\lambda \int_\mathbb{T} \ell_s^a \, \mu(da)} \, d\ell_s^y \right] = \mathbb{E}^y \left[ \int_0^\infty e^{-\lambda \int_\mathbb{T} \ell_s^a \, \mu(da)} \, d\ell_s^x \right]$$

for all $x, y \in \mathbb{R}$. This suffices to establish the result. $\qquad\square$

## 3. Hitting times of bilateral birth-and-death processes

In order to compute certain hitting time distributions for $\xi$, we recall and develop some of the connections between Laplace transforms of hitting times for a birth-and-death process and continued fractions. See Section 12 for some relevant background and notation for continued fractions.

The connection between birth-and-death processes and continued fractions has already been explored, for instance, in [23, 31]. The role of continued fractions in this setting is to pick out the correct solutions of the (generalized) Sturm-Liouville equations [63], whose relationship to quasidiffusions in general is well-laid out in [51].

Suppose that $Z$ is a bilateral birth-and-death process. That is, $Z$ is a continuous time Markov chain on the integers $\mathbb{Z}$ that only makes $\pm 1$ jumps. We assume for concreteness that $Z$ is killed if it reaches $\pm\infty$ in finite time, although this assumption does not feature in the recurrences we derive in this section. Write $\beta_n$ (resp. $\delta_n$) for the rate of jumping to state $n+1$ (resp. $n-1$) from state $n$.

For $n \in \mathbb{Z}$, let $\tau_n = \inf\{t \geq 0 : Z_t = n\}$ be the hitting time of $n$, with the usual convention that the infimum of the empty set is $+\infty$. Set

$$H_n^\downarrow(\lambda) := \mathbb{E}^n[e^{-\lambda \tau_{n-1}}],$$
$$H_n^\uparrow(\lambda) := \mathbb{E}^n[e^{-\lambda \tau_{n+1}}],$$
$$H_{n,m}(\lambda) := \mathbb{E}^n[e^{-\lambda \tau_m}].$$

Note that

$$H_{n,m}(\lambda) = H_n^\uparrow(\lambda) H_{n+1}^\uparrow(\lambda) \cdots H_{m-1}^\uparrow(\lambda), \quad m > n,$$

and

$$H_{n,m}(\lambda) = H_n^\downarrow(\lambda) H_{n-1}^\downarrow(\lambda) \cdots H_{m+1}^\downarrow(\lambda), \quad m < n,$$

and so the fundamental objects to consider are $H_n^\downarrow$ and $H_n^\uparrow$.

Conditioning on the direction of the first jump, we get the recurrence

$$H_n^\downarrow(\lambda) = \mathbb{E}^n \left[ e^{-\lambda \tau_{n-1}} \mathbf{1}_{\tau_{n-1} < \tau_{n+1}} + e^{-\lambda \tau_{n+1}} \mathbf{1}_{\tau_{n+1} < \tau_{n-1}} \mathbb{E}^{n+1} \left[ e^{-\lambda \tau_{n-1}} \right] \right],$$
$$= \frac{\delta_n}{\delta_n + \beta_n + \lambda} + \frac{\beta_n}{\delta_n + \beta_n + \lambda} H_{n+1}^\downarrow(\lambda) H_n^\downarrow(\lambda),$$

which, putting $\rho_n := \frac{\delta_n}{\beta_n}$, can be rearranged as a pair of recurrences

$$(3.1) \qquad H_n^{\downarrow}(\lambda) = \frac{\rho_n}{1 + \rho_n + \frac{\lambda}{\beta_n} - H_{n+1}^{\downarrow}(\lambda)},$$

$$(3.2) \qquad H_{n+1}^{\downarrow}(\lambda) = 1 + \rho_n + \frac{\lambda}{\beta_n} - \frac{\rho_n}{H_n^{\downarrow}(\lambda)}.$$

This leads to two families of terminating continued fractions that connect the Laplace transforms $H_n^{\downarrow}$ for different values of $n$, namely,

$$H_n^{\downarrow}(\lambda) = \cfrac{\rho_n}{1 + \rho_n + \frac{\lambda}{\beta_n} - \cfrac{\rho_{n+1}}{1 + \rho_{n+1} + \frac{\lambda}{\beta_{n+1}} - \cdots - H_{n+m+1}^{\downarrow}}}$$

and

$$H_n^{\downarrow}(\lambda) = 1 + \rho_{n-1} + \cfrac{\lambda}{\beta_{n-1}} - \cfrac{\rho_{n-1}}{1 + \rho_{n-2} + \frac{\lambda}{\beta_{n-2}} - \cdots - \frac{\rho_{n-m-1}}{H_{n-m-1}^{\downarrow}}}.$$

By exchanging $\delta_n$ and $\beta_n$ we get similar relations for $H_n^{\uparrow}$,

$$H_n^{\uparrow}(\lambda) = \frac{1}{1 + \rho_n + \frac{\lambda}{\beta_n} - \rho_n H_{n-1}^{\uparrow}(\lambda)},$$

$$\rho_n H_{n-1}^{\uparrow}(\lambda) = 1 + \rho_n + \frac{\lambda}{\beta_n} - \frac{1}{H_n^{\uparrow}(\lambda)}.$$

If we define

$$s_n(z) := \frac{-\rho_n}{1 + \rho_n + \frac{\lambda}{\beta_n} + z} \quad \text{and} \quad \hat{s}_n(z) := \frac{-\rho_n^{-1}}{1 + \rho_n^{-1} + \frac{\lambda}{\delta_n} + z},$$

then we can write the resulting four continued fraction recurrences as

$$(3.3) \qquad -H_n^{\downarrow}(\lambda) = s_n \circ s_{n+1} \circ \cdots \circ s_{n+m-1}(-H_{n+m}^{\downarrow}(\lambda)),$$

$$(3.4) \qquad -\frac{1}{H_n^{\downarrow}(\lambda)} = \hat{s}_{n-1} \circ \hat{s}_{n-2} \circ \cdots \circ \hat{s}_{n-m}(-\frac{1}{H_{n-m}^{\downarrow}(\lambda)}),$$

$$(3.5) \qquad -\frac{1}{H_n^{\uparrow}(\lambda)} = s_{n+1} \circ s_{n+2} \circ \cdots \circ s_{n+m}(-\frac{1}{H_{n+m}^{\uparrow}(\lambda)}),$$

$$(3.6) \qquad -H_n^{\uparrow}(\lambda) = \hat{s}_n \circ \hat{s}_{n-1} \circ \cdots \circ \hat{s}_{n-m+1}(-H_{n-m}^{\uparrow}(\lambda)).$$

In the context of a unilateral birth-and-death chain (that is, the analogue of our process $Z$ on the state space $\mathbb{N}$), the context considered in [23, 31], there is theory giving conditions under which such continued fractions converge and their classical values give the corresponding Laplace transform. In the bilateral case, not all of the above continued fraction expansions can converge to the classical values, because that would imply, for instance, that $H_n^{\downarrow}(\lambda) = (H_{n-1}^{\uparrow}(\lambda))^{-1}$, but two Laplace transforms of sub-probability measures can only be the reciprocals of each other if both are identically 1, which is certainly not the case here.

In the next section we consider bilateral chains arising from instances of our process $\xi$ on $\mathbb{T}$ and discuss circumstances in which Laplace transforms of hitting times are indeed given by their putative continued fraction representations.

## 4. Hitting times on a scattered subset of $\mathbb{T}$

Suppose in this section that for some $a \in \mathbb{T}$ the infinite set $\mathbb{T} \cap (a, +\infty)$ is discrete with $a$ as an accumulation point. Write $\mathbb{T} \cap (a, +\infty) = \{t_n : n \in \mathbb{Z}\}$ with $t_n < t_{n+1}$ for all $n \in \mathbb{Z}$, and define $\mathcal{Z} : \mathbb{T} \cap (a, +\infty) \to \mathbb{Z}$ by $\mathcal{Z}(t_n) := n$. Then the image under $\mathcal{Z}$ of $\xi$ killed when it exits $(a, +\infty)$ is a bilateral birth-and-death process that can "reach $-\infty$ in finite time and be killed there".

From Proposition 2.2, the jump rates of $Z$ are

$$(4.1) \qquad \delta_n = \frac{1}{(t_n - t_{n-1})(t_{n+1} - t_{n-1})} \text{ and } \beta_n = \frac{1}{(t_{n+1} - t_n)(t_{n+1} - t_{n-1})},$$

$$(4.2) \qquad \text{and so} \quad \rho_n = \frac{(t_{n+1} - t_n)}{(t_n - t_{n-1})}.$$

The convergence properties of the continued fraction expansions given in (3.3)–(3.6) can sometimes be determined by the behavior of $\mathbb{T} \cap (a, +\infty)$ in the neighborhood of its endpoint $a$. We refer the reader to Section 12 for a review of the theory of limit-periodic continued fractions that we use.

It is clear from the construction of $\xi$ as a time change of Brownian motion that $\inf\{t > 0 : \xi_t = a\} = 0$, $\mathbb{P}^a$-a.s. and $\inf\{t > 0 : \xi_t \neq a\} = 0$, $\mathbb{P}^a$-a.s. That is, $a$ is a regular instantaneous point for $\xi$. Thus

$$(4.3) \qquad \lim_{n \to \infty} \mathbb{P}^{-n}[e^{-\lambda \tau_{-n-1}}] = 1$$

and

$$(4.4) \qquad \lim_{n \to \infty} \mathbb{P}^{-n}[e^{-\lambda \tau_{-n+1}} \mid \tau_{-n+1} < \infty] = 1.$$

Note that $\beta_{-n} \to \infty$ as $n \to \infty$. Suppose further that $\rho_{-n} \to \rho \in (1, \infty)$ as $n \to \infty$, Then $\hat{s}_{n-m} \to \hat{s}^*$ as $m \to \infty$, where $\hat{s}^*(z) := \frac{-\rho^{-1}}{1 + \rho^{-1} + z}$ , a transformation with attractive fixed point $-\rho^{-1}$ and repulsive fixed point $-1$. It follows from (4.3) that $\lim_{n \to \infty} H^{\downarrow}_{-n}(\lambda) = 1$. So, by Theorem 12.1, $H^{\downarrow}_{-n}$ is not equal to the classical value of the non-terminating continued fraction corresponding to (3.4). Also, by (4.3),

$$
\begin{aligned}
H^{\uparrow}_{-n}(\lambda) &:= \mathbb{E}^{-n}[e^{-\lambda \tau_{-n+1}}] \\
&= \mathbb{P}^{-n}\{\tau_{-n+1} < \infty\} \mathbb{E}^{-n}[e^{-\lambda \tau_{-n+1}} \mid \tau_{-n+1} < \infty] \\
&= \frac{(t_{-n} - a)}{(t_{-n+1} - a)} \mathbb{E}^{-n}[e^{-\lambda \tau_{-n+1}} \mid \tau_{-n+1} < \infty] \\
&\to \frac{1}{\rho} \quad \text{as} \quad n \to \infty.
\end{aligned}
$$

Thus, Theorem 12.1, applied with indices reversed, implies that the continued fraction expansion in (3.6) converges to the classical value. So, for each $n \in \mathbb{Z}$, $H^{\uparrow}_n = -\lim_{m \to \infty} \hat{s}_n \circ \cdots \circ \hat{s}_{n-m}(0) = \tilde{U}_n / \tilde{U}_{n+1}$, where $\{\tilde{U}_k\}$ is the minimal solution *in the negative direction* to

$$(4.5) \qquad U_{k-1} = (1 + \rho_k^{-1} + \frac{\lambda}{\delta_k}) U_k - \rho_k^{-1} U_{k+1}.$$

In particular, this says that the Laplace transform of the upwards hitting times for the process killed at $a$ are given by a simple formula in terms of the $\{\tilde{U}_m\}$,

$$H_{n,n+m}(\lambda) = \prod_{k=0}^{m-1} H_{n+k}^{\uparrow}(\lambda) = \frac{\tilde{U}_{n+m+1}}{\tilde{U}_{n+1}}, \quad m > 0.$$

Suppose now that $\beta_n$ and $\delta_n$ converge to 0 as $n \to \infty$ in such a way that $\rho_n \to \rho \in (1, \infty)$. An equivalence transformation of continued fractions relate the continued fraction implied by the recurrence (3.3) and the continued fraction implied by the equivalent recurrence

$$(4.6) \qquad -\beta_{n-1} H_n^{\downarrow}(\lambda) = \frac{-\beta_{n-1}\delta_n}{\beta_n + \delta_n + \lambda - \beta_n H_{n+1}^{\downarrow}(\lambda)}.$$

Since $\beta_n$ and $\delta_n$ tend to zero as $n \to \infty$, the limiting transformation is singular and the fixed points tend to zero and $-\lambda$. Since $0 < H_n^{\downarrow}(\lambda) < 1$, $\lim_{n\to\infty} \beta_{n-1} H_n^{\downarrow}(\lambda) = 0$ for all $\lambda > 0$, which is the attractive fixed point of the transformation. Theorem 12.1 implies that the continued fraction converges to the classical value, which is given by the ratio of the minimal solution in the positive direction of the recurrence

$$(4.7) \qquad V_{k+1} = (\beta_k + \delta_k + \lambda) V_k - \beta_{k-1} \delta_k V_{k-1}.$$

However, $U_0 := V_0$ and

$$U_k := \begin{cases} V_k / \prod_{i=0}^{k-1} \beta_i, & k > 0, \\ V_k \prod_{i=k}^{-1} \beta_i, & k < 0, \end{cases}$$

defines a one-one correspondence between solutions to (4.7) and solutions to

$$(4.8) \qquad U_{k+1} = (1 + \rho_k + \frac{\lambda}{\beta_k}) U_k - \rho_k U_{k-1}.$$

Note the sequence $U$ is not the same as in (4.5). Since this correspondence maps minimal solutions to minimal solutions, if we denote by $\{\tilde{U}_k\}$ the minimal solution in the positive direction to (4.8), then

$$H_n^{\downarrow}(\lambda) = \frac{\tilde{U}_n}{\tilde{U}_{n-1}}$$

and

$$H_{n,n-m}(\lambda) = \prod_{k=0}^{m-1} H_{n-k}^{\downarrow}(\lambda) = \frac{\tilde{U}_{n-m-1}}{\tilde{U}_{n-1}}, \quad m > 0.$$

## 5. Introducing the process on $\mathbb{T}_q$

**Note: For the remainder of the paper, we restrict attention to the state space $\mathbb{T} = \mathbb{T}_q := \{q^n \ : \ n \in \mathbb{Z}\} \cup \{-q^n \ : \ n \in \mathbb{Z}\} \cup \{0\}$ for some $q > 1$.**

In this case the measure $\mu$ defining the time change that produces $\xi$ from Brownian motion is given by $\mu = \mu^q$, where $\mu^q(\{q^n\}) = (q^{n+1} - q^{n-1})/2$, $\mu^q(\{-q^n\}) = \mu^q(\{q^n\})$, and $\mu^q(\{0\}) = 0$. Let $\hat{\xi}$ denote the Markov process on $\mathbb{T}_q \cap (0, \infty) = \{q^k :$

$k \in \mathbb{Z}\}$ with distribution starting at $x$ which is that of $\xi$ started at $x$ and killed when it first reaches 0.

By Proposition 2.2 the generator $\mathcal{G}$ of $\xi$ is defined for all $f \in C_0(\mathbb{T}_q)$ for which the following is well-defined and defines a function in $C_0(\mathbb{T}_q)$,

$$(5.1) \qquad (\mathcal{G}f)(x) := \begin{cases} \frac{1}{c_q}\left( \frac{f(qx)}{x^2} + \frac{qf(q^{-1}x)}{x^2} - \frac{(1+q)f(x)}{x^2} \right), & x \in \mathbb{T}_q \setminus \{0\} \\ \lim_{n \to \infty} \frac{1}{2} \frac{f(q^{-n}) + f(-q^{-n}) - 2f(0)}{q^{-2n}}, & x = 0, \end{cases}$$

where $c_q := q^{-1}(q-1)^2(1+q)$. In particular, when our process is at any point $x \neq 0$, it waits for an exponential time with rate proportional to $x^{-2}$ and then jumps further from 0 with probability $1/(1+q)$ or closer to 0 with probability $q/(1+q)$.

We first reinforce our claim that the process $\xi$ on $\mathbb{T}_q$ is a reasonable $q$-analogue of Brownian motion by showing that $\xi$ converges to Brownian motion as the parameter $q$ goes to 1.

**Proposition 5.1.** *For each $q$ let $x_q \in \mathbb{T}_q$ be such that $x_q \to x$ as $q \downarrow 1$. Then the distribution of $\xi$ started at $x_q$ converges as $q \downarrow 1$ (with respect to the usual Skorohod topology on the space of real-valued càdlàg paths) to the distribution of Brownian motion started at $x$.*

*Proof.* Let $(B_t)_{t \in \mathbb{R}_+}$ be a standard Brownian motion with $B_0 = 0$ and let $\ell_t^a$ denote the jointly continuous local time process of $B$.

Set
$$A_u^{\mu^q} := \int_{\mathbb{R}} \ell_u^{a - x_q} \, \mu^q(da)$$
and
$$\theta_t^{\mu^q} := \inf\{u \; : \; A_u^{\mu^q} > t\}.$$

Then the process $(x_q + B(\theta_t^{\mu^q}))_{t \in \mathbb{R}_+}$ has the distribution of $\xi$ under $\mathbb{P}^{x_q}$.

Since $\mu^q$ converges vaguely to the Lebesgue measure $m$ on $\mathbb{R}$ as $q \downarrow 1$, we have
$$\lim_{q \downarrow 1} A_u^{\mu^q} = \int_{\mathbb{R}} \ell_u^{a - x} \, m(da) = u$$
uniformly on compact intervals almost surely, and hence
$$\lim_{q \downarrow 1} \theta_t^{\mu^q} = t$$
uniformly on compact intervals almost surely. Thus $x_q + B(\theta_t^{\mu^q})$ converges to $x + B_t$ uniformly on compact intervals (and hence in the Skorohod topology) almost surely. □

The following lemma shows that $\xi$ obeys a scaling property similar to that of Brownian motion.

**Lemma 5.1.** *The distribution of the process $(\xi_t)_{t \in \mathbb{R}_+}$ under $\mathbb{P}_x$ is the same as that of $(\frac{1}{q}\xi_{q^2t})_{t \in \mathbb{R}_+}$ under $\mathbb{P}_{qx}$. A similar result holds for the killed process $\hat{\xi}$.*

*Proof.* The claim for the process $\xi$ is immediate by checking that properties (I'), (II) and (III) hold for $(\frac{1}{q}\xi_{q^2t})_{t \in \mathbb{R}_+}$. Alternatively, one can verify that the generators of the two processes agree, or use the time-change construction of $\xi$ from Brownian motion and the scaling properties of Brownian motion. The claim for the killed process follows immediately. □

Perhaps the easiest things to calculate about the distribution of $\xi$ are the moments of $\xi_t$. Formally applying the formula for the generator of $\xi$ from Proposition 2.2 to the function $f(x) = x^k$ gives $\mathcal{G}f(x) = \frac{q^{1-k}}{c_q}(1-q^k)(1-q^{k-1})x^{k-2}$. As for the particular cases of $k = 1, 2$ considered in the proof of Proposition 2.1, we can use Dynkin's formula (2.6) and an approximation argument to get the recursion formula

$$\mathbb{E}^x[\xi_t^k] = x^k + \int_0^t \mathbb{E}^x(\mathcal{G}x^k)(\xi_s)\,ds = x^k + \int_0^t \frac{q^{1-k}}{c_q}(1-q^k)(1-q^{k-1})\mathbb{E}^x[\xi_s^{k-2}]\,ds,$$

and hence, using the notation introduced in Section 11,

$$\mathbb{E}^x[\xi_t^k] = \sum_{\substack{m=0 \\ 2|(k-m)}}^{k} c_q^{-\frac{k-m}{2}} \frac{(q;q)_k}{(q;q)_m} q^{\frac{m^2-k^2}{4}} \frac{t^{\frac{k-m}{2}}}{\left(\frac{k-m}{2}\right)!} x^m.$$

where we mean that the sum goes over all $0 \leq m \leq k$ with the same parity as $k$.

The formula shows that, say for $x = 0$, the $k^{\text{th}}$ (even) moments grow like $q^{\frac{k^2}{4}(1+o(1))}t^{\frac{k}{2}}$. This rate of growth is too fast to guarantee that the moments characterize the distribution of $\xi_t$. Note that some well known distributions have moments with this rate of growth, for example, the standard log-normal distribution has $k^{\text{th}}$ moment $e^{\frac{k^2}{2}}$, as does the discrete measure which assigns mass proportional to $e^{\frac{-k^2}{2}}$ at the points $e^k$, $k \in \mathbb{Z}$, [12], 2.3e.

From Proposition 5.1 we would expect informally that the moments of $\xi_t$ should converge to those of a Brownian motion at time $t$ as $q \downarrow 1$. Recall that $c_q = q^{-1}(q-1)^2(1+q)$ and observe that $\lim_{q \downarrow 1}(q-1)^{-\ell}(q;q)_\ell = (-1)^\ell \ell!$. Therefore, if we take $x_q \in \mathbb{T}_q$ with $\lim_{q \downarrow 1} x_q = x \in \mathbb{R}$, then we have

$$\lim_{q \downarrow 1} \mathbb{E}^{x_q}[\xi_t^k] = \sum_{\substack{m=0 \\ 2|(k-m)}}^{k} \frac{k!}{m!(k-m)!} \left(\frac{1}{2}\right)^{\frac{k-m}{2}} \frac{(k-m)!}{\left(\frac{k-m}{2}\right)!} t^{\frac{k-m}{2}} x^m.$$

We recognize the expression on the right hand side as being indeed the $k^{\text{th}}$ moment of a Gaussian random variable with mean $x$ and variance $t$.

## 6. Hitting time distributions for $\mathbb{T}_q$

We once again stress that for the remainder of the paper we are considering the process $\xi$ on the state space $\mathbb{T}_q$.

The general considerations of Section 4 apply to $\mathbb{T}_q \cap (0, \infty)$. In the notation of that section, $t_n = q^n$ for $n \in \mathbb{Z}$. The death and birth rates for the corresponding bilateral birth-and-death process on $\mathbb{Z}$ are, respectively, $\frac{q^{-2n+1}}{c_q}$ and $\frac{q^{-2n}}{c_q}$, where we recall that $c_q = q^{-1}(q-1)^2(1+q)$.

To avoid the constant appearance of factors of $c_q$ in our results, rather than work with $\xi$ and its counterpart $\hat{\xi}$ killed at 0, we will work with the linearly time-changed processes $X = \xi(c_q\cdot)$ and $\hat{X} = \hat{\xi}(c_q\cdot)$. Of course, conclusions for $X$ and $\hat{X}$ can be easily translated into conclusions for $\xi$ and $\hat{\xi}$.

The corresponding bilateral birth-and-death process on $\mathbb{Z}$ has death and birth rates $\delta_n = q^{-2n+1}$ and $\beta_n = q^{-2n}$. In the notation of Section 4, $\rho = \rho_n = q$ and

$$s_n(z) = \frac{-q}{(1+q) + \lambda q^{2n} + z}.$$

Moreover, we have $\tau_n = \inf\{t \in \mathbb{R}_+ : \hat{X}_t = q^n\}$. Note, by the scaling properties in Lemma 5.1, that $H_0^\downarrow(\lambda) = H_n^\downarrow(q^{-2n}\lambda)$, and $H_0^\uparrow(\lambda) = H_n^\uparrow(q^{-2n}\lambda)$. Moreover, recall that

$$H_{n,n-m}(\lambda) = H_n^\downarrow(\lambda)H_{n-1}^\downarrow(\lambda)\cdots H_{n-m+1}^\downarrow(\lambda)$$

and

$$H_{n,n+m}(\lambda) = H_n^\uparrow(\lambda)H_{n+1}^\uparrow(\lambda)\cdots H_{n+m-1}^\uparrow(\lambda),$$

so to compute $H_{n,n-m}(\lambda)$ and $H_{n,n+m}$ it suffices to compute $H_0^\downarrow(\lambda)$ and $H_0^\uparrow(\lambda)$.

From Section 4 we have

$$(6.1) \qquad H_0^\downarrow(\lambda) = \cfrac{q}{1 + q + \lambda - \cfrac{q}{1 + q + \lambda q^2 - \cfrac{q}{1 + q + \lambda q^4 - \ddots}}},$$

$$(6.2) \qquad H_0^\uparrow(\lambda) = \cfrac{1}{1 + q + \lambda - \cfrac{q}{1 + q + \lambda q^{-2} - \cfrac{q}{1 + q + \lambda q^{-4} - \ddots}}}.$$

Closed-form expressions for continued fractions of this form are listed in Ramanujan's "lost" notebook (see the discussion in [5]), and evaluations for various ranges of the parameters (although not all the values we need) can be found in [5, 32, 33] (although in the last several parameter restrictions are omitted).

**Theorem 6.1.** (i) *The Laplace transform of the time to go from 1 to $q^{-1}$ for both $X$ and $\hat{X}$ is*

$$H_0^\downarrow(\lambda) = \frac{q}{\lambda} \frac{{}_0\phi_1(-; 0; q^{-1}; \frac{1}{\lambda q})}{{}_0\phi_1(-; 0; q^{-1}; \frac{1}{\lambda q^{-1}})}.$$

*An alternative expression is*

$$H_0^\downarrow(\lambda) = \frac{1}{(\lambda q^{-1} + 1)} \frac{{}_1\phi_1(0; -\frac{1}{\lambda q}; q^{-2}; -\frac{1}{\lambda q^2})}{{}_1\phi_1(0; -\frac{1}{\lambda q^{-1}}; q^{-2}; -\frac{1}{\lambda q})}.$$

(ii) *The Laplace transform of the time to go from 1 to $q$ for $\hat{X}$ is*

$$H_0^\uparrow(\lambda) = \frac{1}{(q + \lambda)} \frac{{}_1\phi_1(0; -\lambda q^{-3}; q^{-2}; q^{-3})}{{}_1\phi_1(0; -\lambda q^{-1}; q^{-2}; q^{-3})}.$$

*Proof.* (i) Consider the first expression. Since the continued fraction (6.1) converges, by Lemma 12.1 and equation (4.6), the value of $q^{-2n}H_n^\downarrow(\lambda)$ is given by the ratio of consecutive terms of the minimal solution to

$$W_{n+1} = ((1+q)q^{-2n} + \lambda)W_n - q^{-4n+3}W_{n-1}.$$

This recurrence is found in [32] (but with their $q$ as our $q^{-1}$), and the minimal solution is shown to be

$$\tilde{U}_n(\lambda) := q^{-2n(n-1)}\left(\frac{1}{q\lambda}\right)^n {}_0\phi_1(-; 0; q^{-1}; \frac{1}{\lambda q^{2n+1}}).$$

For the second expression, we evaluate (6.1) as follows. Set

$$(6.3) \qquad r_n(\lambda) := {}_1\phi_1(0; -\frac{1}{\lambda}q^{-2n-1}; q^{-2}; -\frac{1}{\lambda}q^{-2n-2})e_{q^{-2}}(-\lambda q^{2n-1}),$$

$$(6.4) \qquad \begin{aligned} h_n(\lambda) &:= -q^{-2n+2}\frac{r_n(\lambda)}{r_{n-1}(\lambda)} \\ &= -\frac{q^{-4n+3}}{(\lambda + q^{-2n+1})}\frac{{}_1\phi_1(0; -\frac{1}{\lambda}q^{-2n-1}; q^{-2}; -\frac{1}{\lambda}q^{-2n-2})}{{}_1\phi_1(0; -\frac{1}{\lambda}q^{-2n+1}; q^{-2}; -\frac{1}{\lambda}q^{-2n})}. \end{aligned}$$

Then, from equation (17) in [5] ,

$$(6.5) \qquad h_n = \frac{-q^{-4n+3}}{(1+q)q^{-2n} + \lambda + h_{n+1}}.$$

This transformation tends to a singular transformation as $n \to \infty$ and the fixed points tend to $x = 0$ and $y = -\lambda$. By (6.4), $q^{4n}h_n \to \frac{q^3}{\lambda}$ as $n \to \infty$, so $h_n \to 0$, and convergence to the classical value holds. However, the continued fraction coming from (6.5) is related by an equivalence transformation to the continued fraction coming from the relation

$$q^{2(n-1)}h_n = \frac{-q}{1 + q + \lambda q^{2n} + q^{2n}h_{n+1}}.$$

This is exactly what is needed to evaluate (6.1), and hence

$$H_n^{\downarrow}(\lambda) = -q^{2(n-1)}h_n(\lambda) = \frac{r_n(\lambda)}{r_{n-1}(\lambda)}.$$

Note that this also shows that $q^{-n(n-1)}r_n(\lambda)$ is a minimal solution in the positive direction to the recurrence

$$U_{n+1}(\lambda) = ((1+q)q^{-2n} + \lambda)U_n(\lambda) - q^{-4n+3}U_{n-1}(\lambda),$$

and is hence equal to $\tilde{U}_n$ above up to a constant multiple.

(ii) Define

$$r'_n(\lambda) = q^{-n}\,{}_1\phi_1(0; -\lambda q^{2n-3}; q^{-2}; q^{-3})/e_{q^{-2}}(-\lambda q^{2n-3}),$$

$$\begin{aligned} g_n(\lambda) &= \frac{r'_{-n+1}}{r'_{-n}} \\ &= -\frac{1}{q}(1 + \lambda q^{-2n-1})\frac{{}_1\phi_1(0; -\lambda q^{-2n-1}; q^{-2}; q^{-3})}{{}_1\phi_1(0; -\lambda q^{-2n-3}; q^{-2}; q^{-3})}. \end{aligned}$$

Equation (13) in [5] simplifies to

$$(6.6) \qquad g_n(\lambda) = \frac{-q}{1 + q + \lambda q^{-2n} + g_{n+1}(\lambda)}.$$

The fixed points of the limiting transformation are $-1$ and $-q$, and

$$\lim_{n \to \infty} g_n(\lambda) = -\frac{1}{q}\frac{{}_1\phi_1(0; 0; q^{-2}; q^{-1})}{{}_1\phi_1(0; 0; q^{-2}; q^{-1})} = -\frac{1}{q}.$$

Hence, by Theorem 12.1, $g_n(\lambda)$ is equal to the classical value of the continued fraction implied by (6.6), and $r'_n$ is a minimal solution in the negative direction to the recursion

$$U_{n+1} = (1 + \frac{1}{q} + \lambda q^{2n-1})U_n - \frac{1}{q}U_{n-1}. \qquad \qquad \square$$

Let $\tau_{-\infty}$ denote the death time of $\hat{X}$. Equivalently, $\tau_{-\infty}$ is the first hitting time of 0 by $X$. Write $H_{n,-\infty}(\lambda) := \mathbb{E}^{q^n}[e^{-\lambda \tau_{-\infty}}]$.

**Corollary 6.1.** *The Laplace transforms of various hitting times for $\hat{X}$ are given by*

$$H_{n,n-m}(\lambda) = \frac{q^{m^2-2mn}}{\lambda^m} \frac{{}_0\phi_1(-;0;q^{-1};\frac{1}{\lambda q^{2n+1}})}{{}_0\phi_1(-;0;q^{-1};\frac{1}{\lambda q^{2(n-m)-1}})}$$

$$= \frac{1}{(-\lambda q^{2n-1};q^{-2})_m} \frac{{}_1\phi_1(0;-\frac{1}{\lambda q^{2n+1}};q^{-2};-\frac{1}{\lambda q^{2n+2}})}{{}_1\phi_1(0;-\frac{1}{\lambda q^{2n-1}}q^{2m};q^{-2};-\frac{1}{\lambda q^{2n}}q^{2m})},$$

$$H_{n,-\infty}(\lambda) = {}_1\phi_1(0;-\frac{1}{\lambda q^{2n+1}};q^{-2};-\frac{1}{\lambda q^{2n+2}})e_{q^{-2}}(-\lambda q^{2n-1})/e_{q^{-2}}(\frac{1}{q}),$$

*and*

$$H_{n,n+m}(\lambda) = \frac{1}{q^m(-\lambda q^{2n+2m-3};q^{-2})_m} \frac{{}_1\phi_1(0;-\lambda q^{2n-3};q^{-2};q^{-3})}{{}_1\phi_1(0;-\lambda q^{2n+2m-3};q^{-2};q^{-3})}.$$

*Proof.* The only result that requires proof is that for $H_{n,-\infty}(\lambda)$. However, by (11.2)

$$\lim_{n\to\infty} {}_1\phi_1(0;-\frac{1}{\lambda}q^{-2n-1};q^{-2};-\frac{1}{\lambda}q^{-2n-2}) = {}_1\phi_0(0;-;q^{-2};\frac{1}{q}) = e_{q^{-2}}(\frac{1}{q}). \qquad \square$$

We can apply known identities to obtain alternatives to the expressions for the Laplace transforms in Theorem 6.1 and Corollary 6.1. For example, equation (13) in [5] gives

$$H_0^\uparrow(\lambda) = 1 - \frac{{}_1\phi_1(0;-q^{-1}\lambda;q^{-2};q^{-1})}{{}_1\phi_1(0;-q^{-1}\lambda;q^{-2};q^{-3})}.$$

Similarly, equation (17) in [5] gives

$$H_0^\downarrow(\lambda) = \frac{q}{q + \lambda \frac{{}_1\phi_1(0;-\frac{1}{q\lambda};q^{-2};-\frac{1}{\lambda})}{{}_1\phi_1(0;-\frac{1}{q\lambda};q^{-2};-\frac{1}{q^2\lambda})}}.$$

The relation

$$(w;q)_\infty \, {}_1\phi_1(0;w;q;c) = (c;q)_\infty \, {}_1\phi_1(0;c;q;w)$$

follows from (III.1) in [30] upon sending $b \to 0$, letting $a = w/z$, and sending $z \to 0$. Similarly, the recurrence

$${}_1\phi_1(0;-\lambda q^{k-4};q^{-2};q^{-3})$$

$$= \frac{(-\frac{1}{\lambda q^{k-2}};q^{-2})_\infty}{(-\frac{1}{\lambda q^{k-1}};q^{-2})_\infty} \, {}_1\phi_1(0;-\lambda q^{k-3};q^{-2};q^{-1})$$

$$- \frac{(-\frac{1}{\lambda q^{k-2}},q^{-1},q^{-1};q^{-2})_\infty}{(-\frac{1}{\lambda q^{k-1}},-\lambda q^{k-2},-\lambda q^{k-3};q^{-2})_\infty} \, {}_1\phi_1(0;-\frac{1}{\lambda q^k};q^{-2};-\frac{1}{\lambda q^{k+1}})$$

comes from (III.31) in [30] by sending $b \to 0$, letting $a = w/z$, and sending $a \to 0$. Both of these identities can be used to obtain alternative formulae for $H_0^\downarrow$ and $H_0^\uparrow$.

We can invert the Laplace transform $H_{0,-\infty}$ in Corollary 6.1 to obtain the distribution of the time $\tau_{-\infty}$ for $X$ to hit 0 starting from $q^n$. Note first of all that

$$\frac{1}{(-\lambda q^{2n-1}; q^{-2})_\infty} = \prod_{i=0}^\infty \frac{1}{1 + \lambda q^{2n-1}q^{-2i}} = \prod_{i=0}^\infty \frac{q^{2i-2n+1}}{q^{2i-2n+1} + \lambda}.$$

Similarly,

$$_1\phi_1(0; -\frac{1}{\lambda q^{2n+1}}; q^{-2}; -\frac{1}{\lambda q^{2n+2}})$$

$$= \sum_{k=0}^\infty \left[ (\lambda q^{2n+2})^k q^{k(k-1)} (-\frac{1}{\lambda q^{2n+1}}; q^{-2})_k (q^{-2}; q^{-2})_k \right]^{-1}$$

$$= \sum_{k=0}^\infty \prod_{l=0}^{k-1} \left[ \frac{q^{-2(l+n)-1}}{(q^{-2(l+n)-1} + \lambda)} \right] \frac{q^{-k}}{(q^{-2}; q^{-2})_k}.$$

Thus under $\mathbb{P}^{q^n}$ the killing time $\tau_{-\infty}$ has the same distribution as the random variable

$$q^{2n} \left( \sum_{i=-\infty}^0 q^{2i-1} T_i + \sum_{i=1}^N q^{2i-1} T_i \right) = q^{2n} \sum_{i=-\infty}^N q^{2i-1} T_i$$

where the $T_i$ are independent rate 1 exponentials and $N$ is distributed according to a $q$-analogue of the Poisson distribution [42], namely,

$$\mathbb{P}\{N = k\} = \frac{1}{e_{q^{-2}}(\frac{1}{q})} \frac{q^{-k}}{(q^{-2}; q^{-2})_k}, \quad k \geq 0.$$

It follows that under $\mathbb{P}^{q^n}$ the distribution of $\tau_{-\infty}$ is also that of the random variable

$$q^{2n+2N-1} \sum_{i=0}^\infty q^{-2i} T_i.$$

A partial fraction expansion of the Laplace transform shows that a convolution of exponential distributions, where the $i^{\text{th}}$ has rate $\alpha_i$, has density

$$t \mapsto \sum_i \alpha_i e^{-\alpha_i t} \prod_{j \neq i} \frac{\alpha_j}{\alpha_j - \alpha_i}.$$

Hence $\sum_{j=0}^\infty q^{-2j} T_j$ has density

$$\begin{aligned}
(6.7) \quad f(t) &:= \sum_{j=0}^\infty q^{2j} e^{-q^{-2j}t} \prod_{k=0}^{j-1} \left( \frac{1}{1 - q^{-2(k-j)}} \right) \prod_{k=j+1}^\infty \left( \frac{1}{1 - q^{-2(k-j)}} \right) \\
&= \frac{1}{(q^{-2}; q^{-2})_\infty} \sum_{j=0}^\infty \frac{q^{2j} e^{-q^{2j}t}}{(q^2; q^2)_j} \\
&= e_{q^{-2}}(-q^{-2}) \sum_{j=0}^\infty \frac{(-1)^j q^{-j(j-1)} e^{-q^{2j}t}}{(q^{-2}; q^{-2})_j}.
\end{aligned}$$

We note in passing that the random variable $\sum_{j=0}^\infty q^{-2j} T_j$ has the same distribution as the exponential functional of the Poisson process $I_{q^{-2}}$ investigated in [4] (see also [3]).

The following result is now immediate.

**Proposition 6.1.** *Under* $\mathbb{P}^{q^n}$, *the hitting time of 0 for X has density*

$$\frac{1}{e_{q^{-2}}(\frac{1}{q})} \sum_{m=0}^{\infty} \frac{q^{-m}}{(q^{-2};q^{-2})_m} q^{2(m+n)+1} f(tq^{2(n+m)+1}), \quad t > 0,$$

*where* $f(t)$ *is defined in (6.7).*

Recall that for Brownian motion started at 1, the hitting time of 0 has the stable($\frac{1}{2}$) density

$$\frac{1}{\sqrt{2\pi t^3}} \exp\left(-\frac{1}{2t}\right), \quad t > 0.$$

It follows from Proposition 5.1 that the distribution of $c_q q^{2N} \sum_{i=0}^{\infty} q^{-2i} T_i$ converges to this stable distribution as $q \downarrow 1$. From Lai's strong law of large numbers for Abelian summation [53] we have that

$$\lim_{q \downarrow 1} \sum_{i=0}^{\infty} (1-q^{-2}) q^{-2i} T_i = \mathbb{E}[T_0] = 1, \quad \text{a.s.}$$

and so $c_q(1-q^{-2})^{-1} q^{2N}$ also converges to the same stable distribution. Taking logarithms, we obtain the following result.

**Proposition 6.2.** *As* $q \downarrow 1$, *the distribution of the random variable*

$$2(\log q)N + \log(q-1)$$

*converges to the distribution with density*

$$\frac{1}{\sqrt{2\pi}} \exp\left(-\frac{1}{2}(x + \exp(-x))\right), \quad -\infty < x < \infty.$$

## 7. Excursion theory for $\mathbb{T}_q$

Recall that $X = \xi(c_q \cdot)$ under $\mathbb{P}^x$ has the same distribution as $(B(\theta_{c_q t}))_{t \in \mathbb{R}_+}$, where $B$ is a Brownian motion started at $x$ with local time process $\ell$, $\theta$ is the right-continuous inverse of the continuous additive functional $A_t = \int \ell_t^a \mu(da)$, and $\mu^q$ is the measure supported on $\mathbb{T}_q$ that is defined by $\mu^q(\{q^n\}) = (q^{n+1} - q^{n-1})/2$, $\mu^q(\{-q^n\}) = \mu(\{q^n\})$, and $\mu^q(\{0\}) = 0$.

Recall also that 0 is a regular instantaneous point for $X$. Thus $X$ has a continuous local time $L$ at 0 that is unique up to constant multiples. We can (and will) take $L_t = \ell_{\theta_{c_q t}}^0$. The inverse of the local time is a subordinator (that is, an increasing Lévy process). Also, there is a corresponding Itô decomposition with respect to the local time of the path of $X$ into a Poisson process of excursions from 0. In this section we determine both the distribution of the subordinator (by giving its Lévy exponent) and the intensity measure of the Poisson process of excursions.

We begin with the following result, which is immediate from Lemma 2.1.

**Lemma 7.1.** *The process X is reversible with respect to the measure* $\mu^q$. *In particular,* $\mu^q$ *is a stationary measure for X.*

We use the excursion theory set-up described in Section VI.8 of [61], which we now briefly review to fix notation. Adjoin an extra cemetery state $\partial$ to $\mathbb{T}_q$. An excursion from 0 is a càdlàg function $f : \mathbb{R}_+ \to \mathbb{T}_q \cup \{\partial\}$ such that $f(0) = 0$ and $f(t) = \partial$ for $t \geq \zeta$, where $\zeta := \inf\{t > 0 : f(t) = \partial$ or $f(t-) = 0\} > 0$. Write $U$ for the space of excursion paths from 0. Using the local time $L$, we can decompose that paths of $X$ under $\mathbb{P}^x$ into a Poisson point process on $\mathbb{R}_+ \times U$ with intensity measure of the form $m \otimes \mathbf{n}$, where $m$ is Lebesgue measure and $\mathbf{n}$ is a $\sigma$-finite measure on $U$ called the Itô excursion measure. The measure $\mathbf{n}$ is time-homogeneous Markov with transition dynamics those of $X$ killed on hitting 0 (and then being sent to $\partial$). Thus $\mathbf{n}$ is completely described by the family of entrance laws $\mathbf{n}_t$, $t > 0$, where

$$\mathbf{n}_t(\Gamma) := \mathbf{n}(\{f \in U : f(t) \in \Gamma\}), \quad \Gamma \subset \mathbb{T}_q.$$

Let $R_\lambda$ denote the $\lambda$-*resolvent* of $X$ for $\lambda > 0$. That is,

$$R_\lambda(x, \Gamma) := \int_0^\infty e^{-\lambda t} \mathbb{P}^x\{X_t \in \Gamma\}\, dt$$

for $x \in \mathbb{T}_q$ and $\Gamma \subseteq \mathbb{T}_q$. In order to identify $\mathbf{n}$, we begin with the following general excursion theory identity (see equation (50.3) in Section VI.8 of [61]).

$$(7.1) \qquad \kappa_\lambda \int_0^\infty e^{-\lambda t} \mathbf{n}_t(\{y\})\, dt = R_\lambda(0, \{y\})$$

where

$$(7.2) \qquad\qquad \kappa_\lambda := \mathbb{E}^0\left[\int_0^\infty e^{-\lambda s}\, dL_s\right].$$

Now, setting $T_x := \inf\{t \in \mathbb{R}_+ : X_t = x\}$, $x \in \mathbb{T}_q$,

$$
\begin{aligned}
R_\lambda(0, \{y\}) &= \lim_{x\to 0} R_\lambda(x, \{y\}) \\
&= \lim_{x\to 0} \frac{\mu^q(\{y\})}{\mu^q(\{x\})} R_\lambda(y, \{x\}) \\
&= \lim_{x\to 0} \frac{\mu^q(\{y\})}{\mu^q(\{x\})} \mathbb{E}^y\left[e^{-\lambda T_x}\right] R_\lambda(x, \{x\}) \\
&= \lim_{x\to 0} \frac{\mu^q(\{y\})}{\mu^q(\{x\})} \mathbb{E}^y\left[e^{-\lambda T_x}\right] \frac{R_\lambda(0, \{x\})}{\mathbb{E}^0\left[e^{-\lambda T_x}\right]} \\
&= \mu^q(\{y\}) \mathbb{E}^y\left[e^{-\lambda T_0}\right] \lim_{x\to 0} \frac{R_\lambda(0, \{x\})}{\mu^q(\{x\})} \\
&= \mu^q(\{y\}) \frac{1}{c_q} \mathbb{E}^y\left[e^{-\lambda T_0}\right] \lim_{x\to 0} \mathbb{E}^0\left[\int_0^\infty e^{-\frac{\lambda}{c_q} A_s}\, d\ell_s^x\right] \\
&= \mu^q(\{y\}) \mathbb{E}^y\left[e^{-\lambda T_0}\right] \frac{1}{c_q} \mathbb{E}^0\left[\int_0^\infty e^{-\frac{\lambda}{c_q} A_s}\, d\ell_s^0\right] \\
&= \mu^q(\{y\}) \mathbb{E}^y\left[e^{-\lambda T_0}\right] \mathbb{E}^0\left[\int_0^\infty e^{-\lambda t}\, d\ell_{\theta_{c_q t}}^0\right],
\end{aligned}
$$

where we used Lemma 7.1 in the second and sixth lines, and a change of variable in the final line.

Thus,

$$\int_0^\infty e^{-\lambda t} \mathbf{n}_t(\{y\})\, dt = \mu^q(\{y\}) \mathbb{E}^y\left[e^{-\lambda T_0}\right],$$

so that

$$\mathbf{n}_t(\{y\}) = \mu^q(\{y\}) \frac{\mathbb{P}^y\{T_0 \in dt\}}{dt}.$$

Now $\mathbb{E}^{q^n}\left[e^{-\lambda T_0}\right] = H_{n,-\infty}(\lambda)$, and so we obtain the following from Corollary 6.1.

**Proposition 7.1.** *The family of entrance laws* $(\mathbf{n}_t)_{t>0}$ *is characterized by*

$$\int_0^\infty e^{-\lambda t} \mathbf{n}_t(\{q^n\}) \, dt$$
$$= \frac{1}{2} q^{n-1}(q^2-1) \, {}_1\phi_1\left(0; -\frac{1}{\lambda q^{2n+1}}; q^{-2}; -\frac{1}{\lambda q^{2n+2}}\right) e_{q^{-2}}(-\lambda q^{2n-1}) \Big/ e_{q^{-2}}\left(\frac{1}{q}\right)$$

*and* $\mathbf{n}_t(\{-q^n\}) = \mathbf{n}_t(\{q^n\})$, $n \in \mathbb{Z}$.

Let $\gamma$ denote the right-continuous inverse of the local time $L$, so that $\gamma$ is a subordinator. Thus $\mathbb{E}^0[e^{-\lambda \gamma_t}] = e^{-t\psi(\lambda)}$ for some Laplace exponent $\psi$.

**Proposition 7.2.** *The distribution of the subordinator* $\gamma$ *is characterized by*

$$\psi(\lambda) = \frac{\lambda (q^2-1) \, (-\frac{1}{\lambda}, -\lambda q^{-2}; q^{-2})_\infty}{q \, (-\frac{\lambda}{q}, -\frac{1}{\lambda q}; q^{-2})_\infty} = \frac{\lambda (q^2-1) \, e_{q^{-2}}(-\frac{\lambda}{q}) e_{q^{-2}}(-\frac{1}{\lambda q})}{q \, e_{q^{-2}}(-\frac{1}{\lambda}) e_{q^{-2}}(-\lambda q^{-2})}.$$

*Proof.* We note the relationship

$$\kappa_\lambda = \mathbb{E}^0\left[\int_0^\infty e^{-\lambda s} \, dL_s\right] = \mathbb{E}^0\left[\int_0^\infty e^{-\lambda \gamma_t} \, dt\right]$$
$$= \int_0^\infty e^{-t\psi(\lambda)} \, dt = \frac{1}{\psi(\lambda)}.$$

Hence, from equation (7.1),

$$\psi(\lambda) = \lambda \int_0^\infty e^{-\lambda t} \mathbf{n}_t(\mathbb{T}_q \setminus \{0\}) \, dt = 2\lambda \int_0^\infty e^{-\lambda t} \mathbf{n}_t(\mathbb{T}_q \cap (0, \infty)) \, dt$$
$$= 2\lambda \sum_{n \in \mathbb{Z}} \frac{1}{2} q^{n-1}(q^2-1) \, {}_1\phi_1\left(0; -\frac{1}{\lambda q^{2n+1}}; q^{-2}; -\frac{1}{\lambda q^{2n+2}}\right)$$
$$\times e_{q^{-2}}(-\lambda q^{2n-1}) \Big/ e_{q^{-2}}\left(\frac{1}{q}\right).$$

Using the following identity to simplify the sum,

$$(-\lambda q^{2n-1}; q^{-2})_\infty = (-\frac{\lambda}{q}; q^{-2})_\infty (-\lambda q^{2n-1}; q^{-2})_n$$
$$= (-\frac{\lambda}{q}; q^{-2})_\infty (-\frac{1}{q\lambda}; q^{-2})_n \, q^{n^2} \lambda^n,$$

we can write part of the above as

$$\sum_{n \in \mathbb{Z}} \frac{q^n \, {}_1\phi_1(0; -\frac{1}{\lambda q^{2n+1}}; q^{-2}; -\frac{1}{\lambda q^{2n+2}})}{(-\lambda q^{2n-1}; q^{-2})_\infty}$$
$$= \frac{1}{(-\frac{\lambda}{q}; q^{-2})_\infty} \sum_{n \in \mathbb{Z}} \frac{q^{-n(n-1)} \lambda^{-n}}{(-\frac{1}{q\lambda}; q^{-2})_n} \, {}_1\phi_1\left(0; -\frac{1}{\lambda q^{2n+1}}; q^{-2}; -\frac{1}{\lambda q^{2n+2}}\right)$$
$$= \frac{1}{(-\frac{\lambda}{q}; q^{-2})_\infty} \sum_{n \in \mathbb{Z}, \, k \geq 0} \frac{q^{-n(n-1)-k(k-1)-k(2n+2)} \lambda^{-n-k}}{(-\frac{1}{q\lambda}; q^{-2})_n (-\frac{q^{-2n}}{\lambda q}; q^{-2})_k (q^{-2}; q^{-2})_k}$$

and by changing indices we get

$$\sum_{n\in\mathbb{Z},\,k\geq0}\frac{q^{-n(n-1)-k(k-1)-k(2n+2)}\lambda^{-n-k}}{(-\frac{1}{q^n};q^{-2})_n(-\frac{q^{-2n}}{\lambda q};q^{-2})_k(q^{-2};q^{-2})_k}$$

$$=\sum_{m\in\mathbb{Z}}\frac{q^{-m(m-1)}\lambda^{-m}}{(-\frac{1}{\lambda q};q^{-2})_m}\sum_{k\geq0}\frac{q^{-2k}}{(q^{-2};q^{-2})_k}$$

$$={}_0\psi_1(-;-\frac{1}{\lambda q};q^{-2};-\frac{1}{\lambda})e_{q^{-2}}(q^{-2}).$$

Recall that $e_{q^{-2}}(z)=1/(z;q^{-2})_\infty$. Moreover, using equation (11.3), we can rewrite the ${}_0\psi_1$ as a product,

$${}_0\psi_1(-;-\frac{1}{\lambda q};q^{-2};-\frac{1}{\lambda})=\frac{(q^{-2},-\frac{1}{\lambda},-\lambda q^{-2};q^{-2})_\infty}{(-\frac{1}{\lambda q},\frac{1}{q};q^{-2})_\infty}.$$

The result now follows.                                                      □

It follows from the scaling property Lemma 5.1 and the uniqueness of the local time at 0 up to a constant multiple that $(L_{q^2t})_{t\in\mathbb{R}_+}$ has the same distribution under $\mathbb{P}^0$ as a constant multiple of $L$. Consequently, the exponent $\psi$ must satisfy the scaling relation $\psi(q^{-2}\lambda)=c\psi(\lambda)$ for some constant $c$. Note from the formula in Proposition 7.2 that, indeed,

$$\psi(q^{-2}\lambda)=q^{-2}\frac{(1+\frac{1}{\lambda}q^2)(1+\lambda q^{-2})^{-1}}{(1+\frac{\lambda}{q})^{-1}(1+\frac{1}{\lambda q}q^2)}\psi(\lambda)=q^{-1}\psi(\lambda).$$

## 8. Resolvent of the killed process on $\mathbb{T}_q\cap(0,\infty)$

Let $\hat{R}_\lambda$ denote the resolvent of the process $\hat{X}$ on $\mathbb{T}_q\cap(0,\infty)$ killed at 0. Recall that $\hat{X}$ goes from $q^n$ to $q^{n-1}$ at rate $q^{-2n+1}$ and from $q^n$ to $q^{n+1}$ at rate $q^{-2n}$. Thus the exit time from $q^n$ is exponentially distributed with rate $q^{-2n+1}+q^{-2n}$, the probability of exiting to $q^{n-1}$ is $\frac{q}{q+1}$, and the probability of exiting to $q^{n+1}$ is $\frac{1}{q+1}$. Moreover, $\mathbb{E}^{q^{n-1}}[e^{-\lambda T_n}]=H_{n-1}^\uparrow(\lambda)$ and $\mathbb{E}^{q^{n+1}}[e^{-\lambda T_n}]=H_{n+1}^\downarrow(\lambda)$. From the strong Markov property we get the recurrence

$$\hat{R}_\lambda(q^n,\{q^n\})$$
$$=\frac{1}{\lambda+q^{-2n+1}+q^{-2n}}+\hat{R}_\lambda(q^n,\{q^n\})$$
$$\times\left(\frac{q^{-2n+1}+q^{-2n}}{\lambda+q^{-2n+1}+q^{-2n}}\frac{q}{q+1}H_{n-1}^\uparrow(\lambda)+\frac{q^{-2n+1}+q^{-2n}}{\lambda+q^{-2n+1}+q^{-2n}}\frac{1}{q+1}H_{n+1}^\downarrow(\lambda)\right),$$

so that

$$\hat{R}_\lambda(q^n,\{q^n\})$$
$$=\left\{\lambda+(q^{-2n+1}+q^{-2n})\left[1-\left(\frac{q}{q+1}H_{n-1}^\uparrow(\lambda)+\frac{1}{q+1}H_{n+1}^\downarrow(\lambda)\right)\right]\right\}^{-1}.$$

Substitute any of the explicit formulae for $H_{n-1}^\uparrow$ and $H_{n+1}^\downarrow$ from Section 6 to get an expression for the on-diagonal terms of the resolvent in terms of basic hypergeometric functions.

To obtain the off-diagonal terms, use the observation

$$\hat{R}_\lambda(q^m, \{q^n\}) = \mathbb{E}^{q^m}\left[e^{-\lambda T_{q^n}}\right]\hat{R}_\lambda(q^n, \{q^n\}) = H_{m,n}(\lambda)\hat{R}_\lambda(q^n, \{q^n\})$$

and then substitute in explicit formulae for $H_{m,n}(\lambda)$ from Section 6 to get expressions in terms of basic hypergeometric functions.

Ideally, one would like to invert the Laplace transform implicit in the resolvent to obtain expressions for the transition probabilities $\mathbb{P}^x\{\hat{X}_t = y\}$. We have not been able to do this.

## 9. Resolvent for $\mathbb{T}_q$

Recall that $R_\lambda$ is the resolvent of the process $X$. The resolvents $R_\lambda$ and $\hat{R}_\lambda$ are related by the equations

$$R_\lambda(x, \{y\}) = \begin{cases} \hat{R}_\lambda(x, \{y\}) + \mathbb{E}^x\left[e^{-\lambda T_0}\right]R_\lambda(0, \{y\}), & x, y \in \mathbb{R}_+, \\ \mathbb{E}^x\left[e^{-\lambda T_0}\right]R_\lambda(0, \{y\}), & x \geq 0, \ y < 0, \\ R_\lambda(-x, \{-y\}), & x < 0. \end{cases}$$

Recall equation (7.1), which says that $R_\lambda(0, \{y\}) = \kappa_\lambda \int_0^\infty e^{-\lambda t}\mathbf{n}_t(\{y\})\, dt$. We know from the proof of Proposition 7.2 that $\kappa_\lambda = \frac{1}{\psi(\lambda)}$ and the statement of Proposition 7.2 gives a simple expression for $\psi(\lambda)$ as a ratio of infinite products. Proposition 7.1 gives an expression for $\int_0^\infty e^{-\lambda t}\mathbf{n}_t(\{y\})\, dt$ in terms of basic hypergeometric functions. Again noting that $\mathbb{E}^{q^m}\left[e^{-\lambda T_{q^n}}\right] = H_{m,n}(\lambda)$, we substitute in explicit formulae for $H_{m,n}(\lambda)$ from Section 6 to get expressions for $R_\lambda(x, \{y\})$ in terms of basic hypergeometric functions.

## 10. A remark on spectral representations

An alternative approach to finding explicit formulae for the quantities of interest would be to find a spectral representation for the generator. This is well-described for general quasidiffusions by Küchler and Salminen [51], who build on the spectral theory of strings [13, 41, 63]. Once one has found solutions to the Sturm-Liouville equation $Gu = -\lambda u$ with appropriate boundary conditions, and the orthogonalizing (spectral) measure, one can write down explicit formulae.

One possible method for carrying this out is to use the well-known spectral representation of transition probabilities of a unilateral birth-and-death process, for which the appropriate eigenfunctions are a family of orthogonal polynomials (see, for example, [37, 38, 39, 62]). If we kill $X$ at $q^{-n}$ for some $n \in \mathbb{Z}$ to obtain a process on $\{q^{-n+1}, q^{-n+2}, \ldots\}$, then the corresponding unilateral birth-and-death process has a specialization of the *associated continuous dual q-Hahn polynomials* as its related family of orthogonal polynomials [32]. However, we have not been able to "take limits as $n \to \infty$" in the resulting spectral representation of the transition probabilities to obtain similar formulae for $\hat{X}$.

Note that our expression for the density of the hitting time to zero of Proposition 6.1 appears to be close to a spectral decomposition – compare to Theorem 3.1 in [51], which gives the density as

$$\frac{1}{\pi}\int_{\mathbb{R}} e^{-\lambda^2 t}C(x; \lambda)\rho(d\lambda)$$

where $C$ is a particular solution to the Sturm-Liouville equation and $\rho$ is the spectral measure. However, to put our expression in this form, the two summations need to be exchanged, which is not straightforward.

## 11. Background on basic hypergeometric functions

For the sake of completeness and to establish notation, we review some of the facts we need about *basic hypergeometric functions* (otherwise known as *q-hypergeometric functions*). For a good tutorial, see the article [45] or the books [2, 30]. In order to make the notation in our review coincide with what is common in the literature, **take $0 < q < 1$ in this section** (this $q$ usually corresponds to $q^{-2}$ in the rest of the paper).

Define the *q-shifted factorial* by

$$(z;q)_n := \prod_{k=0}^{n-1}(1 - zq^k) \quad \text{for } n \in \mathbb{N},\ z \in \mathbb{C},$$

$$(z;q)_\infty := \prod_{k=0}^{\infty}(1 - zq^k) \quad \text{for } |z| < 1.$$

The definition of $(z;q)_n$ may be extended consistently by setting

$$(z;q)_k = \frac{(z;q)_\infty}{(zq^k;q)_\infty} \quad \text{for } k \in \mathbb{Z},\ z \in \mathbb{C}.$$

It will be convenient to use the notation

$$(a_1, a_2, \ldots, a_r; q)_k = (a_1;q)_k(a_2;q)_k \ldots (a_r;q)_k.$$

The *q–hypergeometric series* are indexed by nonnegative integers $r$ and $s$, and for any $\{a_i\} \subset \mathbb{C}$, $\{b_j\} \subset \mathbb{C} \setminus \{q^{-k}\}_{k \geq 0}$ are defined by the series

$$_r\phi_s(a_1, \ldots, a_r; b_1, \ldots, b_s; q; z) := \sum_{k=0}^{\infty} \frac{(a_1, \ldots, a_r; q)_k((-1)^k q^{\frac{k(k-1)}{2}})^{1+s-r} z^k}{(b_1, \ldots, b_s; q; q)_k}.$$

Note the factor $(q;q)_k$ on the bottom, which is not present in the definition used by some authors. The series converges for all $z$ if $r \leq s$, on $|z| < 1$ if $r = s+1$, and only at $z = 0$ if $r > s+1$. Using the property that

$$\lim_{a \to \infty} \frac{(a;q)_n}{a^n} = (-1)^n q^{\frac{n(n-1)}{2}},$$

we get the following useful limit relationships

(11.1) $\lim_{a \to \infty}{}_{r+1}\phi_s(a, a_1, \ldots, a_r; b_1, \ldots, b_s; q; \frac{z}{a}) = {}_r\phi_s(a_1, \ldots, a_r; b_1, \ldots, b_s; q; z),$

(11.2) $\lim_{b \to \infty}{}_r\phi_{s+1}(a_1, \ldots, a_r; b, b_1, \ldots, b_s; q; bz) = {}_r\phi_s(a_1, \ldots, a_r; b_1, \ldots, b_s; q; z),$

as long as the limits stay within the range on which the series converge.

**Theorem 11.1** (The $q$-binomial theorem)**.**

$$_1\phi_0(a; -; q; z) = \frac{(az; q)_\infty}{(z; q)_\infty} \quad \text{if } |z| < 1, |q| < 1, a \in \mathbb{C}.$$

There are (at least) two commonly used $q$–analogues of the exponential function.

$$e_q(z) := \; _1\phi_0(0; -; q; z) = \frac{1}{(z; q)_\infty} = \sum_{k=0}^{\infty} \frac{z^k}{(q; q)_k}, \quad \text{for } |z| < 1,$$

and

$$\begin{aligned} E_q(z) := \; _0\phi_0(-; -; q; -z) &= \frac{1}{e_q(-z)} = (-z; q)_\infty \\ &= \sum_{k=0}^{\infty} \frac{q^{k(k-1)/2}(-z)^k}{(q; q)_k}, \quad \text{for } z \in \mathbb{C}. \end{aligned}$$

The *bilateral $q$-hypergeometric series* also appear in our results. They are defined by

$$_r\psi_s(a_1, \ldots, a_r; b_1, \ldots, b_s; q; z) := \sum_{k=-\infty}^{\infty} \frac{(a_1, \ldots, a_r; q)_k((-1)^k q^{\frac{k(k-1)}{2}})^{s-r} z^k}{(b_1, \ldots, b_s; q)_k}.$$

The sum converges for

$$\begin{cases} \left| \frac{b_1 \cdots b_s}{a_1 \cdots a_r} \right| < |z|, & \text{if } s > r, \\ \left| \frac{b_1 \cdots b_s}{a_1 \cdots a_r} \right| < |z| < 1, & \text{if } s = r, \end{cases}$$

and diverges otherwise.

We use the following extension of the Jacobi triple product identity (see equation (1.49) of [45])

$$(11.3) \quad _0\psi_1(-; c; q; z) := \sum_{k=-\infty}^{\infty} \frac{(-1)^k q^{k(k-1)/2} z^k}{(c; q)_k} = \frac{(q, z, q/z; q)_\infty}{(c, c/z; q)_\infty}, \quad |z| > |c|.$$

## 12. Background on recurrence relations and continued fractions

For nonzero complex numbers $a_n$ and $b_n$, $n \in \mathbb{Z}$, consider the three–term recurrence relation

$$(12.1) \quad U_{n+1} = b_n U_n - a_n U_{n-1}.$$

Its connection to continued fractions can be seen immediately by rearranging to get

$$\frac{U_n}{U_{n-1}} = \frac{a_n}{b_n - \frac{U_{n+1}}{U_n}}.$$

In other words, the sequence $W_n = U_n/U_{n-1}$ solves the recurrence

$$(12.2) \quad W_n W_{n+1} = b_n W_n - a_n.$$

Iterating this recurrence, we get that for any $k \geq 0$,

$$-W_n = \cfrac{-a_n}{b_n - \cfrac{a_{n+1}}{b_{n+1} - \cfrac{a_{n+2}}{b_{n+2} - \cdot\,\cdot\,\cdot\, - \frac{a_{n+k}}{b_{n+k} - W_{n+k}}}}}.$$

We refer to this expression as the *continued fraction expansion associated with the recurrence* (12.2).

A solution $(\tilde{U}_n)_{n \in \mathbb{Z}}$ to (12.1) is said to be a *minimal solution* if, for all linearly independent solutions $V_n$, $\lim_{n \to \infty} |\tilde{U}_n|/|V_n| = 0$. The minimal solution to (12.1), if it exists, is unique up to a constant multiple [56].

For clarity, define the linear fractional transformations

$$s_n(z) = \frac{-a_n}{b_n + z},$$

and write their compositions as $S_m^n = s_{m+1} \circ s_{m+2} \circ \cdots \circ s_n$ and $S^n = S_0^n$. The *classical approximants* to the nonterminating continued fraction (sometimes written $\mathbf{K}[\frac{-a_n}{b_n}]$) are given by

$$\cfrac{-a_1}{b_1 + \cfrac{-a_2}{b_2 + \cfrac{-a_3}{\ddots \, + \frac{-a_n}{b_n}}}} = S^n(0).$$

If we let $P_n$ and $Q_n$ be two solutions to (12.1) with initial conditions $P_{-1} = 1$, $P_0 = 0$, $Q_{-1} = 0$, and $Q_0 = 1$, then it is easy to see that

$$S^n(z) = \frac{P_n + z P_{n-1}}{Q_n + z Q_{n-1}}.$$

The continued fraction $\mathbf{K}[\frac{-a_n}{b_n}]$ is said to *converge* if $S^n(0)$ converges to a (finite) limit as $n$ tends to infinity. If this limit exists, it is called the *classical value* of the continued fraction. However, this is a bit arbitrary, because it can happen, for instance, that for all sequences $(w_n)_{n \in \mathbb{N}}$ that stay away from zero, $S_n(w_n)$ converges to the same limit, different from the limit of $S_n(0)$. The problem is easy to see: suppose that $a_n \to a^*$ and $b_n \to b^*$ as $n \to \infty$, so that $s_n \to s^*$. Each $s_n$ has a pair of fixed points that converge to the fixed points $x$ and $y$ of $s^*$ – suppose $|x| < |y|$, so that $x$ is attractive and $y$ is repulsive. One might imagine that as long as $w_n$ stays away from the repulsive fixed point of $s^*$, then $\lim_{n \to \infty} S_m^n(w_n)$ must converge to $x$ as $m \to \infty$, in which case

$$\lim_{n \to \infty} S_n(w_n) = \lim_{n \to \infty,\, n \geq m} S_m \circ S_m^n(w_n), \quad \forall m \geq 0, \text{ so}$$

$$= \lim_{m \to \infty} \left( \lim_{n \to \infty,\, n \geq m} S_m \circ S_m^n(w_n) \right)$$

$$= \lim_{m \to \infty} S_m(x).$$

Note that just by setting $w_n = S_n^{-1}(z)$, we can get $S_n(w_n)$ converging to any limit in $\mathbb{C}$ we'd like – but to do this, the $w_n$ we choose must converge to the repulsive fixed point. The precise sense in which $w_n$ must "stay away" from the repulsive fixed point is given in the Theorem 12.1 below.

The case in which $a_n \to a^*$ and $b_n \to b^*$, where if $a^* = 0$ then $b^* \neq 0$, is called the *limit 1-periodic* case. Moreover, if the fixed points of $s^*$ are distinct and have different moduli the continued fraction is of *loxodromic type*. All the continued fractions we deal with fall into this category. The following combines Theorem 4 in Chapter II and Theorem 28 in Chapter III of [56]. Here $d(\cdot, \cdot)$ is the spherical metric on $\bar{\mathbb{C}}$.

**Theorem 12.1.** *Let* $K[\frac{-a_n}{b_n}]$ *be limit 1-periodic of loxodromic type.*

(i) *There exists an* $f \in \bar{\mathbb{C}}$ *such that for every sequence* $(w_n)_{n \in \mathbb{N}}$ *for which*

$$\begin{cases} \liminf_{n \to \infty} d(w_n, S_n^{-1}(\infty)) > 0, & \text{when} \quad f \neq \infty, \\ \liminf_{n \to \infty} d(w_n, S_n^{-1}(0)) > 0, & \text{when} \quad f = \infty, \end{cases}$$

*we have* $S_n(w_n) \to f$. *In particular,* $S_n(0) \to f$.

(ii) *Let* $V_0 \in \mathbb{C}$ *and* $V_n = S_n^{-1}(V_0)$ *for* $n > 0$. *Consider* $V_0 \in \mathbb{C}$ *and* $V_n = S_n^{-1}(V_0)$ *for* $n > 0$. *Write* $s^* = \lim_{n \to \infty} s_n$ *and suppose that* $s^*$ *has fixed points* $x, y$ *with* $|x| < |y|$, *so that* $x$ *is attractive and* $y$ *is repulsive for* $s^*$.

- *If* $V_0 = f$, *then* $\lim_{n \to \infty} V_n = x$. *Moreover, if* $f \neq \infty$, *then* $V_0 = \tilde{U}_0 / \tilde{U}_{-1}$, *where* $\tilde{U}_n$ *is a minimal solution to* (12.1).
- *Otherwise,* $\lim_{n \to \infty} V_n = y$.

For a proof, see [56]. This implies the following lemma, which we also use to introduce some more notation. Note that the relation $-W_n = S_n^{-1}(-W_0)$ for $n > 0$ is exactly the relationship implied by (12.2). Note also that this gives explicitly the value of the continued fraction, if it converges to a finite value, in terms of the minimal solution to the associated recurrence relation, a result known as Pincherle's theorem [56].

**Lemma 12.1.** *Suppose that* $(W_n)_{n \in \mathbb{Z}}$ *solves* (12.2) *and that the limits*

$$\beta_{\pm} := \lim_{n \to \infty} \frac{1}{2} \left( b_n \pm \sqrt{b_n^2 - 4a_n} \right),$$

*exist, are finite, and the branches of the square root are chosen so that* $|\beta_-| < |\beta_+|$. *If* $\lim_{n \to \infty} W_n \neq -\beta_+$, *then* $\lim_{n \to \infty} W_n = -\beta_-$, *and for any fixed* $m \in \mathbb{Z}$, *the sequence* $U_n$ *defined by*

$$U_n := \begin{cases} \prod_{k=m+1}^{n} W_k, & n > m, \\ 1, & n = m, \\ \left( \prod_{k=n+1}^{m} W_k \right)^{-1}, & n < m, \end{cases}$$

*is a minimal solution to* (12.1).

*Proof.* Since $\beta_-$ is the limit of the attractive fixed points of the corresponding transformations, and $\beta_+$ is the repulsive fixed point, Theorem 12.1 says that $W_n$ is equal to the classical value of the continued fraction

$$\cfrac{a_n}{b_n - \cfrac{a_{n+1}}{b_{n+1} - \cfrac{a_{n+2}}{b_{n+2} - \ddots}}}$$

if and only if $\lim_{k\to\infty} d(-W_{n+k}, \beta_+) > 0$, so a minimal solution $\tilde{U}_n$ exists, and $W_n = \frac{\tilde{U}_n}{U_{n-1}}$. By definition, $U_n = \frac{\tilde{U}_n}{U_m}$ for all $n \in \mathbb{Z}$. This proves the lemma. □

Two continued fractions are said to be related by an *equivalence transformation* if their sequences of approximants are the same. For example, let $c_k$, $k \in \mathbb{Z}$, be nonzero complex numbers. Since for all $n \geq 0$,

$$\cfrac{a_0}{b_0 + \cfrac{a_1}{b_1 + \cfrac{\cdot}{\cdot \cdot + \cfrac{a_n}{b_n + w_n}}}} = \cfrac{c_0 a_0}{c_0 b_0 + \cfrac{c_0 c_1 a_1}{c_1 b_1 + \cfrac{\cdot}{\cdot \cdot + \cfrac{c_{n-1} c_n a_n}{c_n b_n + c_n w_n}}}},$$

we say that the continued fraction expansions on either side are related by an equivalence transformation.

Note that since we allow the indices in (12.1) and (12.2) to take values in $\mathbb{Z}$, by reversing indices we get another recurrence, another continued fraction, another minimal solution, etc. When we need to distinguish, we will refer to, say, $\tilde{U}_n$ as a minimal solution to (12.1) *in the positive direction* if the above definition holds, and a minimal solution to (12.1) *in the negative direction* if $\lim_{n\to\infty} \tilde{U}_{-n}/V_{-n} = 0$ for some (and hence any) other linearly independent solution $V_n$.

**Acknowledgments.** We thank Jim Pitman for many useful suggestions, and thank Pat Fitzsimmons for suggesting to us that the results of Chacon and Jamison as extended by Walsh could be used to prove Proposition 2.2, thereby strengthening considerably the uniqueness result in an earlier version of the paper.